\documentclass{article}
\usepackage{amssymb}
\usepackage{amsmath}
\usepackage{amscd}

\newtheorem{theorem}{Theorem}

\newtheorem{corollary}[theorem]{Corollary}

\newtheorem{definition}[theorem]{Definition}

\newenvironment{proof}[1][Proof]{\noindent\textbf{#1.} }{\ \rule{0.5em}{0.5em}}

\begin{document}

\title{Further Remarks on Multiple $p$-adic $q$-$L$-Function of Two Variables\footnote{This work was supported by Akdeniz University
Scientific Research Project Unit.}}
\author{Mehmet Cenkci, Yilmaz Simsek\footnote{Corresponding author. e-mail:ysimsek@akdeniz.edu.tr}, Veli Kurt \\
Department of Mathematics, Akdeniz University, 07058-Antalya,
Turkey}
\date{}
\maketitle

\textbf{Abstract : }The object of this paper is to give several properties and applications of the multiple $p$-adic $q$%
-$L$-function of two variables $L_{p,q}^{\left( r\right) }\left(
s,z,\chi\right) $. The explicit formulas relating higher order
$q$-Bernoulli polynomials, which involve sums of products of
higher order $q$-zeta function and higher order Dirichlet
$q$-$L$-function are given. The value of higher order Dirichlet
$p$-adic $q$-$L$-function for positive integers is also
calculated. Furthermore, the Kummer-type congruences for multiple
generalized $q$-Bernoulli polynomials are derived by making use of
the difference theorem of higher order Dirichlet $p$-adic
$q$-$L$-function.

\bigskip


\bigskip


\section{Introduction}

\hspace{0.15in}In \cite{Kim-Cho}, Kim and Cho defined the following multiple
$q$-$L$-function:%
\begin{eqnarray}
L_{q}^{\left( r\right) }\left( s,\chi \right) &=&\frac{1}{%
\prod\limits_{j=1}^{r}\left( s-j\right) }\frac{1}{\left[ F\right] _{q}^{r}}%
\sum_{a_{1},\ldots ,a_{r}=1}^{F}\chi \left( a_{1}+\cdots +a_{r}\right) \left[
a_{1}+\cdots +a_{r}\right] _{q}^{-s+r}  \notag \\
&&\times \sum_{m=0}^{\infty }\binom{r-s}{m}q^{\left( a_{1}+\cdots
+a_{r}\right) m}\left( \frac{\left[ F\right] _{q}}{\left[ a_{1}+\cdots +a_{r}%
\right] _{q}}\right) ^{m}\beta _{m,q^{F}}^{\left( r\right) }.  \label{1.1}
\end{eqnarray}

\noindent They also suggested the following question: ``\textit{Is it
possible to give }$p$\textit{-adic analogue of (\ref{1.1}) which can be
viewed as interpolating, in the same way that }$L_{p,q}\left( s,\chi \right)
$\textit{\ interpolates }$L_{q}\left( s,\chi \right) $\textit{\ in \cite%
{Kim2005}, \cite{Kim2006b}?}''. This question was answered positively by
authors in \cite{Cenkci-Simsek-Kurt} by constructing the following two
variable $p$-adic meromorphic function:%
\begin{eqnarray}
L_{p,q}^{\left( r\right) }\left( s,z,\chi \right) &=&\frac{1}{\left[ F\right]
_{q}^{r}}\frac{1}{\prod\limits_{j=1}^{r}\left( s-j\right) }\sum_{\underset{%
\left( a_{1}+\cdots +a_{r},p\right) =1}{a_{1},\ldots ,a_{r}=1}}^{F}\chi
\left( a_{1}+\cdots +a_{r}\right) \left\langle a_{1}+\cdots +a_{r}+p^{\ast
}z\right\rangle _{q}^{-s+r}  \notag \\
&&\times \sum_{m=0}^{\infty }\binom{r-s}{m}q^{\left( a_{1}+\cdots
+a_{r}+p^{\ast }z\right) m}\left( \frac{\left[ F\right] _{q}}{\left[
a_{1}+\cdots +a_{r}+p^{\ast }z\right] _{q}}\right) ^{m}\beta
_{m,q^{F}}^{\left( r\right) }.  \label{1.2}
\end{eqnarray}

The purpose of this paper is to give further properties of the function $%
L_{p,q}^{\left( r\right) }\left( s,z,\chi \right) $ as well as applications
related to Kummer-type congruences for multiple $q$-Bernoulli polynomials.

Kubota and Leopoldt \cite{Kubota-Leopoldt} proved the existence of
meromorphic function $L_{p}\left( s,\chi \right) $, defined over $p$-adic
number field. $L_{p}\left( s,\chi \right) $ is defined by \cite%
{Ferrero-Greenberg}%
\begin{equation*}
L_{p}\left( s,\chi \right) =\sum_{\underset{\left( n,p\right) =1}{n=1}%
}^{\infty }\frac{\chi \left( n\right) }{n^{s}}=\left( 1-\chi \left( p\right)
p^{-s}\right) L\left( s,\chi \right) ,
\end{equation*}

\noindent where $L\left( s,\chi \right) $ is the Dirichlet $L$-function. $%
L_{p}\left( s,\chi \right) $ interpolates the values%
\begin{equation*}
L_{p}\left( 1-n,\chi \right) =-\frac{1}{n}\left( 1-\chi _{n}\left( p\right)
p^{n-1}\right) B_{n,\chi _{n}}
\end{equation*}

\noindent for $n\in \mathbb{Z}$, $n\geqslant 1$, where $B_{n,\chi }$ denotes
the generalized Bernoulli numbers associated with the primitive Dirichlet
character $\chi $, and $\chi _{n}=\chi \omega ^{-n}$ with $\omega $ being
the Teichm\"{u}ller character (cf. \cite{Diamond}, \cite{Ferrero-Greenberg}, %
\cite{Fox2000}, \cite{Fox2003}, \cite{Iwasawa}, \cite{Kim2002a}, \cite%
{Kim2006a}, \cite{Koblitz1979}, \cite{Koblitz1980}, \cite{Shiratani-Yamamato}%
, \cite{Washington1976}, \cite{Washington1997}).

In \cite{Fox2000,Fox2003}, Fox derived a meromorphic function $L_{p}\left(
s,z,\chi \right) $, which is the two-variable extension of the function $%
L_{p}\left( s,\chi \right) $. Kim \cite{Kim2005} constructed $L_{p,q}\left(
s,z,\chi \right) $, which serves as a $q$-extension of $L_{p}\left( s,z,\chi
\right) $. In \cite{Cenkci-Simsek-Kurt}, the authors defined the multiple $p$%
-adic $q$-$L$-function of two variables $L_{p,q}^{\left( r\right) }\left(
s,z,\chi \right) $, which stands for the higher order generalization of
Kim's $L_{p,q}\left( s,z,\chi \right) $.

Ferrero and Greenberg \cite{Ferrero-Greenberg} evaluated the value $\left(
\partial /\partial s\right) L_{p}\left( 0,\chi \right) $. In \cite{Young},
Young gave an extension of this value by using $p$-adic $L$-function of two
variables $L_{p}\left( s,z,\chi \right) $ under some restrictions on the
character $\chi $. Fox \cite{Fox2003} derived a formula for $\left( \partial
/\partial s\right) L_{p}\left( 0,z,\chi \right) $ without any restrictions.
Kim \cite{Kim2005} evaluated the value $\left( \partial /\partial s\right)
L_{p,q}\left( 0,z,\chi \right) $, which is the $q$-extension and
two-variable extension of the result found by Diamond \cite{Diamond} and
Ferrero and Greenberg \cite{Ferrero-Greenberg}. The authors \cite%
{Cenkci-Simsek-Kurt} obtained a formula for $\left( \partial /\partial
s\right) L_{p,q}^{\left( r\right) }\left( 0,z,\chi \right) $, which
generalizes the results of Kim \cite{Kim2005, Kim2006b}, Fox \cite{Fox2003},
Diamond \cite{Diamond} and Ferrero and Greenberg \cite{Ferrero-Greenberg}.
Further extensions for the value $\left( \partial /\partial s\right)
L_{p}\left( 0,\chi \right) $ can be found in \cite{Kim2002a}, \cite{Kim2006b}%
, \cite{Simsek}.

In recent years, many mathematicians and physicists have investigated zeta
functions, multiple zeta functions, $L$-functions and multiple $q$-Bernoulli
numbers and polynomials because mainly of their interest and importance.
These functions and polynomials are used not only in Complex Analysis and
Mathematical Physics, but also in $p$-adic Analysis and other areas. In
particular, multiple zeta functions and multiple $L$-functions occur within
the context of Knot Theory, Quantum Field Theory, Applied Analysis and
Number Theory (see \cite{Kim2002b}, \cite{Kim2003a}, \cite{Kim2003b}, \cite%
{Kim2006b}, \cite{Nelson-Gartley1994}, \cite{Nelson-Gartley1996}, \cite%
{Srivastava-Kim-Simsek}).

The object of the present sequel to earlier work \cite{Cenkci-Simsek-Kurt}
is to derive several properties and applications of the multiple $p$-adic $q$%
-$L$-function of two variables $L_{p,q}^{\left( r\right) }\left( s,z,\chi
\right) $. We first give a brief summary for zeta and Dirichlet $L$%
-functions and related concepts in preliminary section. In Section
3, we review the definition and construction of the multiple
$q$-$L$-function of two variables and its $p$-adic analogue. We
also find explicit formulas relating higher order $q$-Bernoulli
polynomials, which involve sums of products of $\zeta_{q}^{\left(
r\right) }\left( -n,z_{1}+\cdots+z_{r}\right) $ and $L_{q}^{\left(
r\right) }\left( -n,z_{1}+\cdots+z_{r},\chi \right) $ for integer
$n\geqslant 0$. In Section 4, we evaluate the value
$L_{p,q}^{\left( r\right) }\left( r,z,\chi \right) $ for a positive integer $%
r$ explicitly, from which the value $L_{p,q}\left( 1,z,\chi \right) $ given by Kim \cite%
{Kim2005} is an immediate consequence. In final section, we
purpose to derive Kummer-type congruences for multiple generalized
$q$-Bernoulli polynomials making use of the difference theorem of
$L_{p,q}^{\left( r\right) }\left( s,z,\chi \right) $ and its
consequence, which are already proven in
\cite{Cenkci-Simsek-Kurt}. These congruences are generalizations
of the congruences given by \cite{Carlitz1959},
\cite{Cenkci-Kurt}, \cite{Fox2000}, \cite{Gunaratne1995a,
Gunaratne1995b}, \cite{Jang-Kim-Park}, \cite{Kim2002c},
\cite{Shiratani}.

\section{Preliminaries}

\hspace{0.15in}In complex number field $\mathbb{C}$, the Bernoulli numbers $%
B_{n}$ are defined by means of the generating function%
\begin{equation*}
F\left( t\right) =\frac{t}{e^{t}-1}=\sum_{n=0}^{\infty }B_{n}\frac{t^{n}}{n!}%
,\text{ \ \ \ \ \ \ \ \ \ }\left( \left| t\right| <2\pi \right) .
\end{equation*}%
\noindent $\left\{ B_{n}\right\} $ is the sequence of rational numbers first
considered by Jacob Bernoulli in the study of finite sums of a given power
of consecutive integers. It follows from the generating function definition
that%
\begin{equation*}
B_{0}=1\text{, }B_{1}=-\frac{1}{2}\text{, }B_{2}=\frac{1}{6}\text{, }%
B_{3}=0,B_{4}=-\frac{1}{30},\ldots ,
\end{equation*}%
\noindent and $B_{2k+1}=0$ for $k\in \mathbb{Z}$, $k\geqslant 1$. For an
indeterminate $z$, Bernoulli polynomials $B_{n}\left( z\right) $ are defined
by%
\begin{equation*}
F\left( z,t\right) =\frac{t}{e^{t}-1}e^{zt}=\sum_{n=0}^{\infty }B_{n}\left(
z\right) \frac{t^{n}}{n!},\text{ \ \ \ \ \ \ \ \ \ }\left( \left| t\right|
<2\pi \right) .
\end{equation*}%
\noindent One of the curious facts about Bernoulli numbers and polynomials
is the relation between the Riemann zeta and the Hurwitz (or generalized)
zeta functions.

\begin{theorem}
\label{thm2.1}(\cite{Apostol}) For every integer $n\geqslant 1$,%
\begin{equation*}
\zeta \left( 1-n\right) =-\frac{B_{n}}{n}\text{ and }\zeta \left(
1-n,z\right) =-\frac{B_{n}\left( z\right) }{n},
\end{equation*}%
\noindent where $\zeta \left( s\right) $ and $\zeta \left( s,z\right) $ are
the Riemann and the Hurwitz (or generalized) zeta functions, defined
respectively by%
\begin{equation*}
\zeta \left( s\right) =\sum_{m=1}^{\infty }\frac{1}{m^{s}}\text{ and }\zeta
\left( s,z\right) =\sum_{m=0}^{\infty }\frac{1}{\left( m+z\right) ^{s}},
\end{equation*}%
\noindent with $s\in \mathbb{C}$, $Re\left( s\right) >1$ and $z\in \mathbb{C}
$ with $Re\left( z\right) >0$.
\end{theorem}

For $n\in \mathbb{Z}$, $n\geqslant 1$, a Dirichlet character to the modulus $%
n$ is a multiplicative map $\chi :\mathbb{Z\rightarrow C}$ such that $\chi
\left( a+n\right) =\chi \left( a\right) $ for all $a\in \mathbb{Z}$ and $%
\chi \left( a\right) =0$ if $\left( a,n\right) \neq 1$. Since $a^{\phi
\left( n\right) }\equiv 1\left( \text{mod}n\right) $ for all $a$ such that $%
\left( a,n\right) =1$, $\chi \left( a\right) $ must be the root of unity for
such $a$. If $\chi $ is a Dirichlet character to the modulus $n$, then for
any positive multiple $m$ of $n$, we can induce a Dirichlet character $\psi $
to the modulus $m$ according to $\psi \left( a\right) =\chi \left( a\right) $
if $\left( a,m\right) =1$ and $\psi \left( a\right) =0$ if $\left(
a,m\right) \neq 1$. The minimum modulus $n$ for which a character $\chi $
cannot be induced from some character to the modulus $m$, $m<n$, is called
the conductor of $\chi $, denoted by $f=f_{\chi }$. Throughout, it will be
assumed that each $\chi $ is defined to modulo its conductor. Such a
character is said to be primitive. For primitive Dirichlet characters $\chi $
and $\psi $ having conductors $f_{\chi }$ and $f_{\psi }$, respectively, the
product $\chi \psi $ is defined by $\chi \psi \left( a\right) =\chi \left(
a\right) \psi \left( a\right) $ for all $a\in \mathbb{Z}$ such that $\left(
a,f_{\chi }f_{\psi }\right) =1$. The character $\chi =1$, having conductor $%
f_{1}=1$ is called the principle character.

Among various generalizations of Bernoulli numbers and polynomials,
generalization with a primitive Dirichlet character $\chi $ has a special
case of attention.

\begin{definition}
\label{def1.2}(\cite{Iwasawa}, \cite{Washington1997}) For a primitive
Dirichlet character $\chi $ having conductor $f\in \mathbb{Z}$, $f\geqslant
1 $, the generalized Bernoulli numbers $B_{n,\chi }$ and polynomials $%
B_{n,\chi }\left( z\right) $ associated with $\chi $ are defined by%
\begin{eqnarray*}
F_{\chi }\left( t\right) &=&\sum_{a=1}^{f}\frac{\chi \left( a\right) te^{at}%
}{e^{ft}-1}=\sum_{n=0}^{\infty }B_{n,\chi }\frac{t^{n}}{n!},\text{ \ \ \ \ \
\ \ \ \ }\left( \left| t\right| <\frac{2\pi }{f}\right) \\
F_{\chi }\left( z,t\right) &=&\sum_{a=1}^{f}\frac{\chi \left( a\right)
te^{\left( a+z\right) t}}{e^{ft}-1}=\sum_{n=0}^{\infty }B_{n,\chi }\left(
z\right) \frac{t^{n}}{n!},\text{ \ }\left( \left| t\right| <\frac{2\pi }{f}%
\right)
\end{eqnarray*}%
\noindent respectively.
\end{definition}

\noindent Note that the classical Bernoulli numbers are obtained when $\chi
=1$, in that $B_{n,1}=B_{n}$ if $n\neq 1$ and $B_{1,1}=-B_{1}$. The
generalized Bernoulli numbers and polynomials can be expressed in terms of
Bernoulli polynomials as%
\begin{eqnarray*}
B_{n,\chi } &=&f^{n-1}\sum_{a=1}^{f}\chi \left( a\right) B_{n}\left( \frac{a%
}{f}\right) , \\
B_{n,\chi }\left( z\right) &=&f^{n-1}\sum_{a=1}^{f}\chi \left( a\right)
B_{n}\left( \frac{a+z}{f}\right) .
\end{eqnarray*}

Given a primitive Dirichlet character $\chi $, having conductor $f$, the
Dirichlet $L$-function associated with $\chi $ is defined by (\cite{Apostol}%
, \cite{Washington1997})%
\begin{equation*}
L\left( s,\chi \right) =\sum_{m=1}^{\infty }\frac{\chi \left( m\right) }{%
m^{s}},
\end{equation*}%
\noindent for $s\in \mathbb{C}$, $Re\left( s\right) >1$. It is well known %
\cite{Washington1997} that $L\left( s,\chi \right) $ may be analytically
continued to the whole complex plane, except for a simple pole at $s=1$ when
$\chi =1$, in which case the Riemann zeta function, $\zeta \left( s\right)
=L\left( s,1\right) $ is obtained. The generalized Bernoulli numbers share a
particular relationship with the Dirichlet $L$-function, in that%
\begin{equation*}
L\left( 1-n,\chi \right) =-\frac{B_{n,\chi }}{n}
\end{equation*}%
\noindent for $n\in \mathbb{Z}$, $n\geqslant 1$.

For $r\in \mathbb{Z}$, $r\geqslant 1$, the Bernoulli numbers $B_{n}^{\left(
r\right) }$ and polynomials $B_{n}^{\left( r\right) }\left( z\right) $ of
order $r$ (also called the multiple Bernoulli numbers and polynomials,
respectively) may be defined by means of \cite[Chapter 6]{Norlund}%
\begin{eqnarray*}
F^{\left( r\right) }\left( t\right) &=&\left( \frac{t}{e^{t}-1}\right)
^{r}=\sum_{n=0}^{\infty }B_{n}^{\left( r\right) }\frac{t^{n}}{n!},\ \ \ \ \
\ \ \ \ \left( \left| t\right| <2\pi \right) \\
F^{\left( r\right) }\left( z,t\right) &=&\left( \frac{t}{e^{t}-1}\right)
^{r}e^{zt}=\sum_{n=0}^{\infty }B_{n}^{\left( r\right) }\left( z\right) \frac{%
t^{n}}{n!},\text{ \ \ \ }\left( \left| t\right| <2\pi \right)
\end{eqnarray*}%
\noindent respectively. For $r=1$, the classical Bernoulli numbers and
polynomials are obtained.

Let $\chi $ be a Dirichlet character of conductor $f$. In \cite%
{Jang-Kim-Park}, the multiple generalized Bernoulli numbers $B_{n,\chi
}^{\left( r\right) }$ attached to $\chi $ are defined by%
\begin{equation*}
F_{\chi }^{\left( r\right) }\left( t\right) =\sum_{a_{1},\ldots ,a_{r}=1}^{f}%
\frac{\chi \left( a_{1}+\cdots +a_{r}\right) t^{r}e^{\left( a_{1}+\cdots
+a_{r}\right) t}}{\left( e^{ft}-1\right) ^{r}}=\sum_{n=0}^{\infty }B_{n,\chi
}^{\left( r\right) }\frac{t^{n}}{n!},\ \ \ \ \ \ \ \ \ \left( \left|
t\right| <\frac{2\pi }{f}\right)
\end{equation*}%
\noindent for $r\in \mathbb{Z}$, $r\geqslant 1$. For an indeterminate $z$,
the multiple generalized Bernoulli polynomials $B_{n,\chi }^{\left( r\right)
}\left( z\right) $ attached to $\chi $ are naturally given by%
\begin{equation*}
F_{\chi }^{\left( r\right) }\left( z,t\right) =\sum_{a_{1},\ldots
,a_{r}=1}^{f}\frac{\chi \left( a_{1}+\cdots +a_{r}\right) t^{r}e^{\left(
z+a_{1}+\cdots +a_{r}\right) t}}{\left( e^{ft}-1\right) ^{r}}%
=\sum_{n=0}^{\infty }B_{n,\chi }^{\left( r\right) }\left( z\right) \frac{%
t^{n}}{n!},\ \ \ \ \ \ \ \ \ \left( \left| t\right| <\frac{2\pi }{f}\right) .
\end{equation*}%
\noindent It can be readily seen that%
\begin{equation*}
B_{n,\chi }^{\left( r\right) }=f^{n-r}\sum_{a_{1},\ldots ,a_{r}=1}^{f}\chi
\left( a_{1}+\cdots +a_{r}\right) B_{n}^{\left( r\right) }\left( \frac{%
a_{1}+\cdots +a_{r}}{f}\right)
\end{equation*}%
\noindent and%
\begin{equation*}
B_{n,\chi }^{\left( r\right) }\left( z\right) =f^{n-r}\sum_{a_{1},\ldots
,a_{r}=1}^{f}\chi \left( a_{1}+\cdots +a_{r}\right) B_{n}^{\left( r\right)
}\left( \frac{z+a_{1}+\cdots +a_{r}}{f}\right) .
\end{equation*}

Throughout this paper, $p$ will denote a prime number, $\mathbb{Z}_{p}$, $%
\mathbb{Q}_{p}$, $\overline{\mathbb{Q}}_{p}$ and $\mathbb{C}_{p}$ will be
used to represent, respectively, the $p$-adic integers, the $p$-adic
numbers, the algebraic closure of $\mathbb{Q}_{p}$ and the completion of $%
\overline{\mathbb{Q}}_{p}$ with respect to $p$-adic absolute value $\left|
\cdot \right| _{p}$, which is normalized so that $\left| p\right|
_{p}=p^{-1} $. On $\mathbb{C}_{p}$, the absolute value is non-Archimedean,
and so for any $a$, $b\in \mathbb{C}_{p}$, $\left| a+b\right| _{p}\leqslant
\max \left\{ \left| a\right| _{p},\left| b\right| _{p}\right\} $. We denote
a particular subring of $\mathbb{C}_{p}$ as%
\begin{equation*}
R=\left\{ a\in \mathbb{C}_{p}:\left| a\right| _{p}\leqslant 1\right\} .
\end{equation*}%
\noindent If $z\in \mathbb{C}_{p}$ such that $\left| z\right| _{p}\leqslant
\left| p\right| _{p}^{m}$, where $m\in \mathbb{Q}$, then $z\in p^{m}R$, and
this can be also written as $z\equiv 0\left( mod p^{m}R\right) $.

Let $p^{\ast }=4$ if $p=2$ and $p^{\ast }=p$ otherwise. Note that there
exists $\phi \left( p^{\ast }\right) $ distinct solutions, modulo $p^{\ast }$%
, to the equation $x^{\phi \left( p^{\ast }\right) }-1=0$, and each solution
must be congruent to one of the values $a\in \mathbb{Z}$, where $1\leqslant
a\leqslant p^{\ast }$, $\left( a,p\right) =1$. Thus, given $a\in \mathbb{Z}$
with $\left( a,p\right) =1$, there exists a unique $\omega \left( a\right)
\in \mathbb{Z}_{p}$, where $\omega \left( a\right) ^{\phi \left( p^{\ast
}\right) }=1$, such that $\omega \left( a\right) \equiv a\left( \text{mod}%
p^{\ast }\mathbb{Z}_{p}\right) $. Letting $\omega \left( a\right) =0$ for $%
a\in \mathbb{Z}$ such that $\left( a,p\right) \neq 1$, it can be seen that $%
\omega $ is actually a Dirichlet character having conductor $f_{\omega
}=p^{\ast }$, called the Teichm\"{u}ller character. Let $\left\langle
a\right\rangle =\omega ^{-1}\left( a\right) a$. Then $\left\langle
a\right\rangle \equiv 1\left( \text{mod}p^{\ast }\mathbb{Z}_{p}\right) $
(see also \cite{Iwasawa}, \cite{Kim2005}, \cite{Washington1997}, \cite{Young}%
).

For the context in the sequel, an extension of the definition of the Teichm%
\"{u}ller character is needed. If $z\in \mathbb{C}_{p}$ such that $\left|
z\right| _{p}\leqslant 1$, then for any $a\in \mathbb{Z}$, $a+p^{\ast
}z\equiv a\left( modp^{\ast }R\right) $. Thus, for $z\in \mathbb{C}_{p}$, $%
\left| z\right| _{p}\leqslant 1$, $\omega \left( a+p^{\ast }z\right) =\omega
\left( a\right) $. Also, for these values of $z$, let $\left\langle
a+p^{\ast }z\right\rangle =\omega ^{-1}\left( a\right) \left( a+p^{\ast
}z\right) $ (see \cite{Fox2000}, \cite{Fox2003}).

$q$-extensions of Bernoulli numbers are first studied by Carlitz \cite%
{Carlitz1948}. The corresponding numbers are called as $q$-Bernoulli
numbers, denoted by $\beta _{n,q}$, and defined by means of the symbolic
formula%
\begin{equation*}
\beta _{0,q}=\frac{q-1}{\text{log}q}\text{, \ \ \ \ }\left( q\beta
_{q}+1\right) ^{n}-\beta _{n,q}=\delta _{n,1},
\end{equation*}%
\noindent with the usual convention about replacing $\beta _{q}^{j}$ by $%
\beta _{j,q}$ and $\delta _{n,1}$ is the Kronecker symbol (cf. \cite%
{Cenkci-Can-Kurt}, \cite{Cenkci-Can}, \cite{Kim2002a}, \cite{Kim2002b}, \cite%
{Kim2004a}, \cite{Simsek-Kim-Rim}). Note that, $\underset{q\rightarrow 1}{%
\lim }\beta _{n,q}=B_{n}$, the classical Bernoulli numbers. $q$-Bernoulli
polynomials $\beta _{n,q}\left( z\right) $ are defined by (cf. \cite%
{Kim2002b}, \cite{Kim2005})%
\begin{equation*}
\beta _{n,q}\left( z\right) =\left( q^{z}\beta _{q}+\left[ z\right]
_{q}\right) ^{n}=\sum_{k=0}^{n}\binom{n}{k}q^{kz}\beta _{k,q}\left[ z\right]
_{q}^{n-k},
\end{equation*}%
\noindent where%
\begin{equation*}
\left[ z\right] _{q}=\frac{1-q^{z}}{1-q}.
\end{equation*}%
\noindent When talking about $q$-extensions, $q$ can be variously considered
as an indeterminate, a complex number $q\in \mathbb{C}$ or a $p$-adic number
$q\in \mathbb{C}_{p}$. If $q\in \mathbb{C}$, we assume that $\left| q\right|
<1$. If $q\in \mathbb{C}_{p}$, then $\left| 1-q\right| _{p}<\left| p\right|
_{p}^{1/\left( p-1\right) }=p^{-1/\left( p-1\right) }$. Thus, for $\left|
x\right| _{p}\leqslant 1$, we have $q^{x}=\exp \left( x\log _{p}q\right) $,
where $\log _{p}$ is the Iwasawa $p$-adic logarithm function (see \cite%
{Iwasawa}).

For $q\in \mathbb{C}$, $\left| q\right| <1$, Kim \cite{Kim2005} gave
generating functions of $q$-Bernoulli numbers and polynomials respectively by%
\begin{eqnarray*}
F_{q}\left( t\right) &=&\frac{q-1}{\text{log}q}e^{t/\left( 1-q\right)
}-t\sum_{n=0}^{\infty }q^{n}e^{\left[ n\right] _{q}t}=\sum_{n=0}^{\infty
}\beta _{n,q}\frac{t^{n}}{n!}, \\
F_{q}\left( z,t\right) &=&\frac{q-1}{\text{log}q}e^{t/\left( 1-q\right)
}-t\sum_{n=0}^{\infty }q^{n+z}e^{\left[ n+z\right] _{q}t}=\sum_{n=0}^{\infty
}\beta _{n,q}\left( z\right) \frac{t^{n}}{n!}
\end{eqnarray*}%
\noindent for $\left| t\right| <1$. For a Dirichlet character $\chi $ with
conductor $f$, the generalized $q$-Bernoulli numbers and polynomials
associated with $\chi $ are defined by the rules%
\begin{eqnarray*}
F_{q,\chi }\left( t\right) &=&-t\sum_{n=1}^{\infty }\chi \left( n\right)
q^{n}e^{\left[ n\right] _{q}t}=\sum_{n=0}^{\infty }\beta _{n,q,\chi }\frac{%
t^{n}}{n!}, \\
F_{q,\chi }\left( z,t\right) &=&-t\sum_{n=1}^{\infty }\chi \left( n\right)
q^{n+z}e^{\left[ n+z\right] _{q}t}=\sum_{n=0}^{\infty }\beta _{n,q,\chi
}\left( z\right) \frac{t^{n}}{n!}
\end{eqnarray*}%
\noindent for $\left| t\right| <1$, respectively (cf. \cite{Kim2005}). From
these formulas, it can be easily obtained that%
\begin{equation*}
\beta _{n,q,\chi }\left( z\right) =\sum_{k=0}^{n}\binom{n}{k}q^{kz}\beta
_{k,q,\chi }\left[ z\right] _{q}^{n-k}.
\end{equation*}%
\noindent The generalized $q$-Bernoulli numbers and polynomials can be
expressed in terms of $q$-Bernoulli polynomials as%
\begin{eqnarray*}
\beta _{n,q,\chi } &=&\left[ f\right] _{q}^{n-1}\sum_{a=1}^{f}\chi \left(
a\right) \beta _{n,q^{f}}\left( \frac{a}{f}\right) , \\
\beta _{n,q,\chi }\left( z\right) &=&\left[ f\right] _{q}^{n-1}%
\sum_{a=1}^{f}\chi \left( a\right) \beta _{n,q^{f}}\left( \frac{a+z}{f}%
\right) .
\end{eqnarray*}

\section{Multiple $q$-$L$-Function of Two Variables and its $p$-adic Analogue%
}

\hspace{0.15in}This section is devoted to recall the multiple $q$-$L$%
-function of two variables $L_{q}^{\left( r\right) }\left( s,z,\chi \right) $%
, and its $p$-adic analogue $L_{p,q}^{\left( r\right) }\left( s,z,\chi
\right) $ which were introduced in \cite{Cenkci-Simsek-Kurt}. We also
summarize the analytic continuation, special values and explicit formulas
for these functions.

For $r\in \mathbb{Z}$, $r\geqslant 1$, the multiple $q$-Bernoulli numbers $%
\beta _{n,q}^{\left( r\right) }$ and the multiple $q$-Bernoulli polynomials $%
\beta _{n,q}^{\left( r\right) }\left( z\right) $, are defined respectively
by means of the generating functions (cf. \cite{Kim2004b})%
\begin{equation*}
F_{q}^{\left( r\right) }\left( t\right) =\left( -t\right)
^{r}\sum_{n_{1},\ldots ,n_{r}=0}^{\infty }q^{n_{1}+\cdots +n_{r}}e^{\left[
n_{1}+\cdots +n_{r}\right] _{q}t}=\sum_{n=0}^{\infty }\beta _{n,q}^{\left(
r\right) }\frac{t^{n}}{n!},
\end{equation*}

\noindent and

\begin{equation}
F_{q}^{\left( r\right) }\left( z,t\right) =\left( -t\right)
^{r}\sum_{n_{1},\ldots ,n_{r}=0}^{\infty }q^{z+n_{1}+\cdots +n_{r}}e^{\left[
z+n_{1}+\cdots +n_{r}\right] _{q}t}=\sum_{n=0}^{\infty }\beta _{n,q}^{\left(
r\right) }\left( z\right) \frac{t^{n}}{n!}  \label{3.1}
\end{equation}%
\noindent for $\left| t\right| <1$. The complex analytic multiple $q$-zeta
function was defined in \cite{Kim2004b} as follows:

\begin{definition}
\label{def3.1}For $s\in \mathbb{C}$, $Re\left( s\right) >r$ and $z\in
\mathbb{C}$, $Re\left( z\right) >0$, the multiple $q$-zeta function $\zeta
_{q}^{\left( r\right) }\left( s,z\right) $ is defined by means of%
\begin{equation*}
\zeta _{q}^{\left( r\right) }\left( s,z\right) =\sum_{m_{1},\ldots
,m_{r}=0}^{\infty }\frac{q^{z+m_{1}+\cdots +m_{r}}}{\left[ z+m_{1}+\cdots
+m_{r}\right] _{q}^{s}}.
\end{equation*}
\end{definition}

\noindent For $s\in \mathbb{C}$, the following integral representation was
obtained by Kim \cite{Kim2004b}:%
\begin{eqnarray*}
\zeta _{q}^{\left( r\right) }\left( s,z\right) &=&\frac{1}{\Gamma \left(
s\right) }\int\limits_{0}^{\infty }t^{s-1-r}F_{q}^{\left( r\right) }\left(
z,-t\right) dt \\
&=&\sum_{m_{1},\ldots ,m_{r}=0}^{\infty }q^{z+m_{1}+\cdots +m_{r}}\frac{1}{%
\Gamma \left( s\right) }\int\limits_{0}^{\infty }t^{s-1}e^{-\left[
z+m_{1}+\cdots +m_{r}\right] _{q}t}dt.
\end{eqnarray*}

\noindent This integral representation yields the following assertion:

\begin{theorem}
\label{thm3.2}(\cite{Kim2004b}) For $n\in \mathbb{Z}$, $n\geqslant 0$ and $%
r\in \mathbb{Z}$, $r\geqslant 1$,%
\begin{eqnarray*}
\zeta _{q}^{\left( r\right) }\left( -n,z\right) &=&\left( -1\right) ^{r}%
\frac{n!}{\left( n+r\right) !}\beta _{n+r,q}^{\left( r\right) }\left(
z\right) \\
&=&\frac{\left( -1\right) ^{r}}{\left( n+1\right) _{r}}\beta
_{n+r,q}^{\left( r\right) }\left( z\right) ,
\end{eqnarray*}%
\noindent where $\left( n\right) _{r}$ is the Pochhammer symbol defined by%
\begin{equation*}
\left( n\right) _{0}=1\text{ and }\left( n\right) _{r}=n\left( n+1\right)
\cdots \left( n+r-1\right) .
\end{equation*}
\end{theorem}

For a Dirichlet character $\chi $, the multiple generalized $q$-Bernoulli
numbers $\beta _{n,q,\chi }^{\left( r\right) }$ are defined by (cf. \cite%
{Kim2004b}, \cite{Kim-Cho})%
\begin{eqnarray*}
\hspace{-0.1in}F_{q,\chi }^{\left( r\right) }\left( t\right) &=&\left(
-t\right) ^{r}\sum_{a_{1},\ldots ,a_{r}=1}^{f}\chi \left( a_{1}+\cdots
+a_{r}\right) \sum_{k_{1},\ldots ,k_{r}=0}^{\infty }q^{a_{1}+k_{1}f+\cdots
+a_{r}+k_{r}f}e^{\left[ a_{1}+k_{1}f+\cdots +a_{r}+k_{r}f\right] _{q}t} \\
&=&\left( -t\right) ^{r}\sum_{n_{1},\ldots ,n_{r}=1}^{\infty }\chi \left(
n_{1}+\cdots +n_{r}\right) q^{n_{1}+\cdots +n_{r}}e^{\left[ n_{1}+\cdots
+n_{r}\right] _{q}t} \\
&=&\sum_{n=0}^{\infty }\beta _{n,q,\chi }^{\left( r\right) }\frac{t^{n}}{n!}
\end{eqnarray*}%
\noindent For an indeterminate $z$, this definition can be expanded in order
to define the multiple generalized $q$-Bernoulli polynomials $\beta
_{n,q,\chi }^{\left( r\right) }\left( z\right) $ as follows (cf. \cite%
{Cenkci-Simsek-Kurt}):%
\begin{eqnarray}
&&\hspace{-1in}F_{q,\chi }^{\left( r\right) }\left( z,t\right) =\left(
-t\right) ^{r}\sum_{a_{1},\ldots ,a_{r}=1}^{f}\chi \left( a_{1}+\cdots
+a_{r}\right)  \notag \\
&\times &\sum_{k_{1},\ldots ,k_{r}=0}^{\infty }q^{z+a_{1}+k_{1}f+\cdots
+a_{r}+k_{r}f}e^{\left[ z+a_{1}+k_{1}f+\cdots +a_{r}+k_{r}f\right] _{q}t}
\notag \\
&=&\left( -t\right) ^{r}\sum_{n_{1},\ldots ,n_{r}=1}^{\infty }\chi \left(
n_{1}+\cdots +n_{r}\right) q^{z+n_{1}+\cdots +n_{r}}e^{\left[ z+n_{1}+\cdots
+n_{r}\right] _{q}t}  \notag \\
&=&\sum_{n=0}^{\infty }\beta _{n,q,\chi }^{\left( r\right) }\left( z\right)
\frac{t^{n}}{n!}.  \label{3.2}
\end{eqnarray}%
\noindent From (\ref{3.1}) and (\ref{3.2}), it readily seen that%
\begin{equation}
\beta _{n,q,\chi }^{\left( r\right) }\left( z\right) =\left[ f\right]
_{q}^{n-r}\sum_{a_{1},\ldots ,a_{r}=1}^{f}\chi \left( a_{1}+\cdots
+a_{r}\right) \beta _{n,q^{f}}^{\left( r\right) }\left( \frac{z+a_{1}+\cdots
+a_{r}}{f}\right) .  \label{3.3}
\end{equation}%
\noindent The multiple Dirichlet $q$-$L$-function of two variables is given
by the following definition:

\begin{definition}
\label{def3.3}(\cite{Cenkci-Simsek-Kurt}) For $s\in \mathbb{C}$, $Re\left(
s\right) >r$, $z\in \mathbb{C}$, $Re\left( z\right) >0$ and a Dirichlet
character $\chi $, the multiple Dirichlet $q$-$L$-function of two variables $%
L_{q}^{\left( r\right) }\left( s,z,\chi \right) $ is defined by means of%
\begin{equation*}
L_{q}^{\left( r\right) }\left( s,z,\chi \right) =\sum_{m_{1},\ldots
,m_{r}=0}^{\infty }\frac{\chi \left( m_{1}+\cdots +m_{r}\right)
q^{z+m_{1}+\cdots +m_{r}}}{\left[ z+m_{1}+\cdots +m_{r}\right] _{q}^{s}}.
\end{equation*}
\end{definition}

\noindent Special values of $L_{q}^{\left( r\right) }\left( s,z,\chi \right)
$ are given by the following theorem:

\begin{theorem}
\label{thm3.4}For $n\in \mathbb{Z}$, $n\geqslant 0$ and $r\in \mathbb{Z}$, $%
r\geqslant 1$, we have%
\begin{equation*}
L_{q}^{\left( r\right) }\left( -n,z,\chi \right) =\frac{\left( -1\right) ^{r}%
}{\left( n+1\right) _{r}}\beta _{n+r,q,\chi }^{\left( r\right) }\left(
z\right) .
\end{equation*}
\end{theorem}

\begin{proof}
The proof of this theorem was given in \cite{Cenkci-Simsek-Kurt} by using
complex integration. Here, we give another proof.

Writing $m_{i}=a_{i}+k_{i}f$ with $a_{i}=1,2,\ldots ,f$, $k_{i}=0,1,2,\ldots
$, where $i=1,\ldots ,r$,%
\begin{eqnarray*}
L_{q}^{\left( r\right) }\left( -n,z,\chi \right) &=&\sum_{m_{1},\ldots
,m_{r}=0}^{\infty }\frac{\chi \left( m_{1}+\cdots +m_{r}\right)
q^{z+m_{1}+\cdots +m_{r}}}{\left[ z+m_{1}+\cdots +m_{r}\right] _{q}^{-n}} \\
&=&\sum_{a_{1},\ldots ,a_{r}=1}^{f}\sum_{k_{1},\ldots ,k_{r}=0}^{\infty
}\chi \left( a_{1}+k_{1}f+\cdots +a_{r}+k_{r}f\right)
q^{z+a_{1}+k_{1}f+\cdots +a_{r}+k_{r}f} \\
&&\times \left[ z+a_{1}+k_{1}f+\cdots +a_{r}+k_{r}f\right] _{q}^{n} \\
&=&\sum_{a_{1},\ldots ,a_{r}=1}^{f}\chi \left( a_{1}+\cdots +a_{r}\right)
\sum_{k_{1},\ldots ,k_{r}=0}^{\infty }\left( q^{f}\right) ^{\frac{%
^{z+a_{1}+\cdots +a_{r}}}{f}+k_{1}+\cdots +k_{r}} \\
&&\times \left[ \frac{z+a_{1}+\cdots +a_{r}}{f}+k_{1}+\cdots +k_{r}\right]
_{q^{f}}^{n}\left[ f\right] _{q}^{n}.
\end{eqnarray*}

\noindent By Definition \ref{def3.1}, last expression equals%
\begin{equation*}
L_{q}^{\left( r\right) }\left( -n,z,\chi \right) =\left[ f\right]
_{q}^{n}\sum_{a_{1},\ldots ,a_{r}=1}^{f}\chi \left( a_{1}+\cdots
+a_{r}\right) \zeta _{q^{f}}^{\left( r\right) }\left( -n,\frac{%
z+a_{1}+\cdots +a_{r}}{f}\right) .
\end{equation*}

\noindent Making use of Theorem \ref{thm3.2}, we get%
\begin{equation*}
L_{q}^{\left( r\right) }\left( -n,z,\chi \right) =\left[ f\right] _{q}^{n}%
\frac{\left( -1\right) ^{r}}{\left( n+1\right) _{r}}\sum_{a_{1},\ldots
,a_{r}=1}^{f}\chi \left( a_{1}+\cdots +a_{r}\right) \beta
_{n+r,q^{f}}^{\left( r\right) }\left( \frac{z+a_{1}+\cdots +a_{r}}{f}\right)
.
\end{equation*}

\noindent Now, utilizing equation (\ref{3.3}) (with $n\longmapsto n+r$) we
obtain%
\begin{equation*}
L_{q}^{\left( r\right) }\left( -n,z,\chi \right) =\frac{\left( -1\right) ^{r}%
}{\left( n+1\right) _{r}}\beta _{n+r,q,\chi }^{\left( r\right) }\left(
z\right) ,
\end{equation*}

\noindent the desired result.
\end{proof}

In \cite{Kim2001, Kim2003a}, Kim gave closed expressions for the
sums of products of $q$-Bernoulli numbers and generalized
$q$-Bernoulli numbers. Using these sums, we have nested sums over
combinations of multiple $q$-zeta function and multiple
$q$-Dirichlet $L$-function as follows:

\begin{theorem}
\label{thm3.6}We have%
\begin{eqnarray*}
&&\hspace{-0.3in}\zeta _{q}^{\left( r\right) }\left( -n,z_{1}+\cdots +z_{r}\right)=\frac{\left( -1\right) ^{r}}{\left( n+1\right) _{r}}\sum_{\underset{%
i_{1}+\cdots +i_{r}=n+r}{i_{1},\ldots ,i_{r}\geqslant 0}%
}\\
&&\hspace{-0.3in}\times\sum_{k_{1}=0}^{n+r-i_{1}}\sum_{k_{2}=0}^{n+r-i_{1}-i_{2}}\cdots
\sum_{k_{r-1}=0}^{n+r-i_{1}-i_{2}-\cdots
-i_{r-1}}\binom{n+r}{i_{1},\ldots
,i_{r}}\binom{n+r-i_{1}}{k_{1}}\cdots \binom{n+r-i_{1}-i_{2}-\cdots -i_{r-1}%
}{k_{r-1}} \\
&&\hspace{-0.3in}\times \beta _{k_{1}+i_{1},q}\left( z_{1}\right)
\beta _{k_{2}+i_{2},q}\left( z_{2}\right) \cdots \beta
_{k_{r-1}+i_{r-1},q}\left( z_{r-1}\right) \beta _{i_{r},q}\left(
z_{r}\right) \left( q-1\right) ^{k_{1}+\cdots +k_{r-1}},
\end{eqnarray*}

\noindent and%
\begin{eqnarray*}
&&\hspace{-0.3in}L_{q}^{\left( r\right) }\left( -n,z_{1}+\cdots
+z_{r},\chi\right)=\frac{\left( -1\right) ^{r}}{\left( n+1\right)
_{r}}\left[ f\right] _{q}^{n}\sum_{a_{1},\ldots ,a_{r}=1}^{f}\chi
\left( a_{1}+\cdots +a_{r}\right) \sum_{\underset{i_{1}+\cdots
+i_{r}=n+r}{i_{1},\ldots ,i_{r}\geqslant 0}}\\
&&\hspace{-0.3in}\times\sum_{k_{1}=0}^{n+r-i_{1}}%
\sum_{k_{2}=0}^{n+r-i_{1}-i_{2}}\cdots
\sum_{k_{r-1}=0}^{n+r-i_{1}-i_{2}-\cdots -i_{r-1}}\binom{n+r}{i_{1},\ldots
,i_{r}}\binom{n+r-i_{1}}{k_{1}}\cdots \binom{n+r-i_{1}-i_{2}-\cdots -i_{r-1}%
}{k_{r-1}} \\
&&\hspace{-0.3in}\times \beta _{k_{1}+i_{1},q^{f}}\left(
\frac{a_{1}+z_{1}}{f}\right) \beta _{k_{2}+i_{2},q^{f}}\left(
\frac{a_{2}+z_{2}}{f}\right) \cdots \beta
_{k_{r-1}+i_{r-1},q^{f}}\left( \frac{a_{r-1}+z_{r-1}}{f}\right)\\
&&\hspace{-0.3in}\times\beta _{i_{r},q^{f}}\left(
\frac{a_{r}+z_{r}}{f}\right) \left( q^{f}-1\right) ^{k_{1}+\cdots
+k_{r-1}}.
\end{eqnarray*}
\end{theorem}

\begin{proof}By using the formula of Kim \cite{Kim2001}%
\begin{eqnarray}
&&\hspace{-0.7in}\beta _{n,q}^{\left( r\right) }\left(
z_{1}+\cdots +z_{r}\right)=\sum_{\underset{i_{1}+\cdots +i_{r}=n+r}{i_{1},\ldots ,i_{r}\geqslant 0}%
}\sum_{k_{1}=0}^{n+r-i_{1}}\sum_{k_{2}=0}^{n+r-i_{1}-i_{2}}\cdots
\sum_{k_{r-1}=0}^{n+r-i_{1}-i_{2}-\cdots
-i_{r-1}}\notag \\
&&\hspace{-0.3in}\times \binom{n+r}{i_{1},\ldots
,i_{r}}\binom{n+r-i_{1}}{k_{1}}\cdots \binom{n+r-i_{1}-i_{2}-\cdots -i_{r-1}%
}{k_{r-1}}  \label{3.7} \\
&&\hspace{-0.3in}\times \beta _{k_{1}+i_{1},q}\left( z_{1}\right)
\beta _{k_{2}+i_{2},q}\left( z_{2}\right) \cdots \beta
_{k_{r-1}+i_{r-1},q}\left( z_{r-1}\right) \beta _{i_{r},q}\left(
z_{r}\right) \left( q-1\right) ^{k_{1}+\cdots +k_{r-1}},  \notag
\end{eqnarray}

\noindent and Theorem \ref{thm3.2}, we obtain the first formula. The second
formula follows from (\ref{3.7}), (\ref{3.3}) and Theorem \ref{thm3.4}.
\end{proof}

Let $F\in \mathbb{Z}$, $F\geqslant 1,$ and

\begin{equation*}
H_{q}^{\left( r\right) }\left( s,z:a_{1},\ldots ,a_{r}|F\right) =\sum_{%
\underset{m_{j}\equiv a_{j}\left( \text{mod}F\right) }{m_{1},\ldots ,m_{r}>0}%
}^{\infty }\frac{q^{z+m_{1}+\cdots +m_{r}}}{\left[ z+m_{1}+\cdots +m_{r}%
\right] _{q}^{s}}
\end{equation*}

\noindent be the multiple partial $q$-zeta function (cf. \cite%
{Cenkci-Simsek-Kurt}). For $m_{j}=a_{j}+n_{j}F$, $j=1,\ldots ,r$, it can be
written that%
\begin{eqnarray*}
H_{q}^{\left( r\right) }\left( s,z:a_{1},\ldots ,a_{r}|F\right)
&=&\sum_{n_{1},\ldots ,n_{r}=0}^{\infty }\frac{q^{z+a_{1}+n_{1}F+\cdots
+a_{r}+n_{r}F}}{\left[ z+a_{1}+n_{1}F+\cdots +a_{r}+n_{r}F\right] _{q}^{s}}
\\
&=&\sum_{n_{1},\ldots ,n_{r}=0}^{\infty }\frac{\left( q^{F}\right) ^{\frac{%
z+a_{1}+\cdots +a_{r}}{F}+n_{1}+\cdots +n_{r}}}{\left[ F\right] _{q}^{s}%
\left[ \frac{z+a_{1}+\cdots +a_{r}}{F}+n_{1}+\cdots +n_{r}\right]
_{q^{F}}^{s}} \\
&=&\left[ F\right] _{q}^{-s}\zeta _{q^{F}}^{\left( r\right) }\left( s,\frac{%
z+a_{1}+\cdots +a_{r}}{F}\right) .
\end{eqnarray*}

\noindent Thus, the function $H_{q}^{\left( r\right) }\left(
s,z:a_{1},\ldots ,a_{r}|F\right) $ is a meromorphic function for $s\in
\mathbb{C}$ with poles at $s=1,2,\ldots ,r$ and for $s=-n$, $n\in \mathbb{Z}$%
, $n\geqslant 0$, it is obvious that%
\begin{eqnarray}
H_{q}^{\left( r\right) }\left( -n,z:a_{1},\ldots ,a_{r}|F\right) &=&\left[ F%
\right] _{q}^{n}\zeta _{q^{F}}^{\left( r\right) }\left( -n,\frac{%
z+a_{1}+\cdots +a_{r}}{F}\right)  \notag \\
&=&\left[ F\right] _{q}^{n}\frac{\left( -1\right) ^{r}}{\left( n+1\right)
_{r}}\beta _{n+r,q^{F}}^{\left( r\right) }\left( \frac{z+a_{1}+\cdots +a_{r}%
}{F}\right) .  \label{3.4}
\end{eqnarray}

Let $F$ be a positive integer multiple of the conductor $f$. The multiple
Dirichlet $q$-$L$-function of two variables $L_{q}^{\left( r\right) }\left(
s,z,\chi \right) $ is given in terms of $H_{q}^{\left( r\right) }\left(
s,z:a_{1},\ldots ,a_{r}|F\right) $ as (cf. \cite{Cenkci-Simsek-Kurt})%
\begin{equation*}
L_{q}^{\left( r\right) }\left( s,z,\chi \right) =\sum_{a_{1},\ldots
,a_{r}=1}^{F}\chi \left( a_{1}+\cdots +a_{r}\right) H_{q}^{\left( r\right)
}\left( s,z:a_{1},\ldots ,a_{r}|F\right) .
\end{equation*}

\noindent For $n\in \mathbb{Z}$, $n\geqslant 0$,%
\begin{equation*}
L_{q}^{\left( r\right) }\left( -n,z,\chi \right) =\sum_{a_{1},\ldots
,a_{r}=1}^{F}\chi \left( a_{1}+\cdots +a_{r}\right) H_{q}^{\left( r\right)
}\left( -n,z:a_{1},\ldots ,a_{r}|F\right) ,
\end{equation*}

\noindent thus from (\ref{3.4}), it can be written that%
\begin{equation}
L_{q}^{\left( r\right) }\left( -n,z,\chi \right) =\frac{\left( -1\right) ^{r}%
}{\left( n+1\right) _{r}}\left[ F\right] _{q}^{n}\sum_{a_{1},\ldots
,a_{r}=1}^{F}\chi \left( a_{1}+\cdots +a_{r}\right) \beta
_{n+r,q^{F}}^{\left( r\right) }\left( \frac{z+a_{1}+\cdots +a_{r}}{F}\right)
.  \label{3.5}
\end{equation}

\noindent By substituting the expansion%
\begin{equation*}
\beta _{n+r,q^{F}}^{\left( r\right) }\left( \frac{z+a_{1}+\cdots +a_{r}}{F}%
\right) =\sum_{k=0}^{n+r}\binom{n+r}{k}\beta _{k,q^{F}}^{\left( r\right)
}q^{k\left( z+a_{1}+\cdots +a_{r}\right) }\left[ \frac{z+a_{1}+\cdots +a_{r}%
}{F}\right] _{q^{F}}^{n+r-k},
\end{equation*}

\noindent of multiple $q$-Bernoulli polynomials on the right of (\ref{3.5}),
the equation%
\begin{eqnarray*}
L_{q}^{\left( r\right) }\left( -n,z,\chi \right) &=&\frac{\left( -1\right)
^{r}}{\left( n+1\right) _{r}}\left[ F\right] _{q}^{-r}\sum_{a_{1},\ldots
,a_{r}=1}^{F}\chi \left( a_{1}+\cdots +a_{r}\right) \left[ z+a_{1}+\cdots
+a_{r}\right] _{q}^{n+r} \\
&&\times \sum_{k=0}^{n+r}\binom{n+r}{k}\beta _{k,q^{F}}^{\left( r\right)
}q^{k\left( z+a_{1}+\cdots +a_{r}\right) }\left( \frac{\left[ F\right] _{q}}{%
\left[ z+a_{1}+\cdots +a_{r}\right] _{q}}\right) ^{k}
\end{eqnarray*}

\noindent can be obtained. Last equation yields%
\begin{eqnarray*}
L_{q}^{\left( r\right) }\left( s,z,\chi \right) &=&\frac{1}{\left[ F\right]
_{q}^{r}}\frac{1}{\prod\limits_{j=1}^{r}\left( s-j\right) }%
\sum_{a_{1},\ldots ,a_{r}=1}^{F}\chi \left( a_{1}+\cdots +a_{r}\right) \left[
z+a_{1}+\cdots +a_{r}\right] _{q}^{-s+r} \\
&&\times \sum_{k=0}^{\infty }\binom{-s+r}{k}\beta _{k,q^{F}}^{\left(
r\right) }q^{k\left( z+a_{1}+\cdots +a_{r}\right) }\left( \frac{\left[ F%
\right] _{q}}{\left[ z+a_{1}+\cdots +a_{r}\right] _{q}}\right) ^{k},
\end{eqnarray*}

\noindent which can be used to construct the $p$-adic analogue of the
function $L_{q}^{\left( r\right) }\left( s,z,\chi \right) $ (cf. \cite%
{Cenkci-Simsek-Kurt}).

Let $q\in \mathbb{C}_{p}$ with $\left| 1-q\right| _{p}<p^{-1/\left(
p-1\right) }$. For an arbitrary character $\chi $ and $n\in \mathbb{Z}$, let
$\chi _{n}=\chi \omega ^{-n}$ in the sense of product of characters, where $%
\omega $ being the Teichm\"{u}ller character. Also, let $\left\langle
a\right\rangle _{q}=\omega ^{-1}\left( a\right) \left[ a\right] _{q}$.\ Then
$\left\langle a\right\rangle _{q}\equiv 1\left(modp^{\ast }R\right) $. If $%
z\in \mathbb{C}_{p}$, $\left| z\right| _{p}\leqslant 1$, then for any $a\in
\mathbb{Z}$, it can be written that $\left\langle a+p^{\ast }z\right\rangle
_{q}=\omega ^{-1}\left( a\right) \left[ a+p^{\ast }z\right] _{q}$ so that $%
\left\langle a+p^{\ast }z\right\rangle _{q}\equiv 1\left( modp^{\ast
}R\right) $ for $z\in \mathbb{C}_{p}$, $\left| z\right| _{p}\leqslant 1$
(cf. \cite{Cenkci-Can}, \cite{Cenkci-Simsek-Kurt}, \cite{Kim2005}, \cite%
{Kim2006a}, \cite{Simsek}).

Let%
\begin{equation*}
A_{j}\left( x\right) =\sum_{n=0}^{\infty }a_{n,j}x^{n},
\end{equation*}

\noindent $a_{n,j}\in \mathbb{C}_{p}$, $j=0,1,\ldots ,$ be a sequence of
formal power series, each of which converges in a fixed subset%
\begin{equation*}
D=\left\{ s\in \mathbb{C}_{p}:\left| s\right| _{p}\leqslant \left| p^{\ast
}\right| _{p}^{-1}p^{-1/\left( p-1\right) }\right\}
\end{equation*}

\noindent of $\mathbb{C}_{p}$ such that $a_{n,j}\rightarrow a_{n,0}$ as $%
j\rightarrow \infty $ for all $n$, and for each $s\in D$ and $\epsilon >0$,
there exists $n_{0}=n_{0}\left( s,\epsilon \right) $ such that%
\begin{equation*}
\left| \sum_{n\geqslant n_{0}}^{\infty }a_{n,j}s^{n}\right| _{p}<\epsilon
\end{equation*}

\noindent for all $j$. Then, $\underset{j\rightarrow \infty }{\lim }%
A_{j}\left( s\right) =A_{0}\left( s\right) $ for any $s\in D$. This fact was
used by Washington \cite{Washington1976} to construct the $p$-adic
equivalent of Dirichlet $L$-function $L\left( s,\chi \right) $.

By making use of this method, the authors \cite{Cenkci-Simsek-Kurt}
constructed the multiple $p$-adic $q$-$L$-function of two variables as
follows:

\begin{eqnarray}
L_{p,q}^{\left( r\right) }\left( s,z,\chi \right) &=&\frac{1}{\left[ F\right]
_{q}^{r}}\frac{1}{\prod\limits_{j=1}^{r}\left( s-j\right) }\sum_{\underset{%
\left( a_{1}+\cdots +a_{r},p\right) =1}{a_{1},\ldots ,a_{r}=1}}^{F}\chi
\left( a_{1}+\cdots +a_{r}\right) \left\langle a_{1}+\cdots +a_{r}+p^{\ast
}z\right\rangle _{q}^{-s+r}  \notag \\
&&\times \sum_{k=0}^{\infty }\binom{-s+r}{k}\beta _{k,q^{F}}^{\left(
r\right) }q^{k\left( a_{1}+\cdots +a_{r}+p^{\ast }z\right) }\left( \frac{%
\left[ F\right] _{q}}{\left[ a_{1}+\cdots +a_{r}+p^{\ast }z\right] _{q}}%
\right) ^{k},  \label{3.6}
\end{eqnarray}

\noindent where $F$ is a positive integer multiple of $p^{\ast }$ and the
conductor $f$ of $\chi $. The analytical properties and analytic
continuation of $L_{p,q}^{\left( r\right) }\left( s,z,\chi \right) $ are
given by the following theorem (\cite{Cenkci-Simsek-Kurt}):

\begin{theorem}
\label{thm3.5}Let $F$ be a positive integer multiple of $f$ and $p^{\ast }$,
and $L_{p,q}^{\left( r\right) }\left( s,z,\chi \right) $ be defined as (\ref%
{3.6}). Then, $L_{p,q}^{\left( r\right) }\left( s,z,\chi \right) $ is
analytic for $z\in \mathbb{C}_{p}$, $\left| z\right| _{p}\leqslant 1$,
provided $s\in D$, except $s=1,2,\ldots ,r$ when $\chi =1$, and meromorphic
with simple poles at $s=1,2,\ldots ,r$ when $\chi =1$. Furthermore, for each
$n\in \mathbb{Z}$, $n\geqslant 0$,%
\begin{equation*}
L_{p,q}^{\left( r\right) }\left( -n,z,\chi \right) =\frac{\left( -1\right)
^{r}}{\left( n+1\right) _{r}}\left( \beta _{n+r,q,\chi _{n+r}}^{\left(
r\right) }\left( p^{\ast }z\right) -\chi _{n+r}\left( p\right) \left[ p%
\right] _{q}^{n}\beta _{n+r,q^{p},\chi _{n+r}}^{\left( r\right) }\left(
p^{-1}p^{\ast }z\right) \right) .
\end{equation*}
\end{theorem}

In next sections, we will consider some properties and applications of the
function $L_{p,q}^{\left( r\right) }\left( s,z,\chi \right) $.

\section{The Value of $L_{p,q}^{\left( r\right) }\left( s,z,\protect\chi %
\right) $ at $s=r$}

In this section, we evaluate the value $L_{p,q}^{\left( r\right) }\left(
r,z,\chi \right) $ for a positive integer $r$.

\begin{theorem}
\label{thm4.1}Let $\chi $ be a Dirichlet character of conductor $f$ and let $%
F$ be a positive integer multiple of $f$ and $p^{\ast }$. Then for $r\in
\mathbb{Z}$, $r\geqslant 1$,%
\begin{eqnarray*}
L_{p,q}^{\left( r\right) }\left( r,z,\chi \right)
&=&\frac{1}{\left[ F\right] _{q}^{r}}\frac{1}{\left( r-1\right)
!}\sum_{\underset{\left( a_{1}+\cdots +a_{r},p\right)
=1}{a_{1},\ldots ,a_{r}=1}}^{F}\chi \left( a_{1}+\cdots
+a_{r}\right) \left\{ -\beta _{0,q^{F}}^{\left( r\right) }\log
_{p}\left\langle a_{1}+\cdots
+a_{r}+p^{\ast }z\right\rangle _{q}\right. \\
&&+\left. \sum_{k=r}^{\infty }\frac{\left( -1\right) ^{k}}{k}\beta
_{k,q^{F}}^{\left( r\right) }q^{k\left( a_{1}+\cdots
+a_{r}+p^{\ast }z\right) }\left( \frac{\left[ F\right]
_{q}}{\left[ a_{1}+\cdots +a_{r}+p^{\ast }z\right] _{q}}\right)
^{k}\right\} .
\end{eqnarray*}
\end{theorem}

\begin{proof}
It follows from (\ref{3.6}) that%
\begin{eqnarray*}
&&\hspace{-0.5in}L_{p,q}^{\left( r\right) }\left( s,z,\chi \right) =\frac{1}{%
\left[ F\right] _{q}^{r}}\sum_{\underset{\left( a_{1}+\cdots +a_{r},p\right)
=1}{a_{1},\ldots ,a_{r}=1}}^{F}\chi \left( a_{1}+\cdots +a_{r}\right) \frac{%
\left\langle a_{1}+\cdots +a_{r}+p^{\ast }z\right\rangle _{q}^{-s+r}}{%
\prod\limits_{j=1}^{r}\left( s-j\right) } \\
&&\times \sum_{k=0}^{r-1}\binom{-s+r}{k}\beta _{k,q^{F}}^{\left( r\right)
}q^{k\left( a_{1}+\cdots +a_{r}+p^{\ast }z\right) }\left( \frac{\left[ F%
\right] _{q}}{\left[ a_{1}+\cdots +a_{r}+p^{\ast }z\right] _{q}}\right) ^{k}
\\
&&+\frac{1}{\left[ F\right] _{q}^{r}}\sum_{\underset{\left( a_{1}+\cdots
+a_{r},p\right) =1}{a_{1},\ldots ,a_{r}=1}}^{F}\chi \left( a_{1}+\cdots
+a_{r}\right) \left\langle a_{1}+\cdots +a_{r}+p^{\ast }z\right\rangle
_{q}^{-s+r} \\
&&\times \sum_{k=r}^{\infty }\frac{\binom{-s+r}{k}}{\prod\limits_{j=1}^{r}%
\left( s-j\right) }\beta _{k,q^{F}}^{\left( r\right) }q^{k\left(
a_{1}+\cdots +a_{r}+p^{\ast }z\right) }\left( \frac{\left[ F\right] _{q}}{%
\left[ a_{1}+\cdots +a_{r}+p^{\ast }z\right] _{q}}\right) ^{k}.
\end{eqnarray*}

\noindent Now, by using Taylor expansion at $s=r$, we have%
\begin{eqnarray*}
&&\hspace{-0.5in}\frac{\left\langle a_{1}+\cdots +a_{r}+p^{\ast
}z\right\rangle _{q}^{-s+r}}{\prod\limits_{j=1}^{r}\left( s-j\right) } \\
&&=-\frac{\log _{p}\left\langle a_{1}+\cdots +a_{r}+p^{\ast }z\right\rangle
_{q}}{\prod\limits_{j=1}^{r-1}\left( s-j\right) }\left( 1+\left\{ \text{%
terms involving the powers of }\left( -s+r\right) \right\} \right) ,
\end{eqnarray*}

\noindent so that%
\begin{equation*}
\underset{s\rightarrow r}{\lim }\frac{\left\langle a_{1}+\cdots
+a_{r}+p^{\ast }z\right\rangle _{q}^{-s+r}}{\prod\limits_{j=1}^{r}\left(
s-j\right) }=-\frac{\log _{p}\left\langle a_{1}+\cdots +a_{r}+p^{\ast
}z\right\rangle _{q}}{\left( r-1\right) !},
\end{equation*}

\noindent and%
\begin{equation*}
\hspace{-0.5in}\frac{\binom{-s+r}{k}}{\prod\limits_{j=1}^{r}\left(
s-j\right) }=\frac{\left( -1\right) ^{k}s\left( s+1\right) \cdots \left(
s-r+k-1\right) }{k!},
\end{equation*}

\noindent so that%
\begin{equation*}
\underset{s\rightarrow r}{\lim }\frac{\binom{-s+r}{k}}{\prod%
\limits_{j=1}^{r}\left( s-j\right) }=\frac{1}{\left( r-1\right) !}\frac{%
\left( -1\right) ^{k}}{k}.
\end{equation*}

\noindent We therefore get%
\begin{eqnarray*}
&&\hspace{-0.7in}L_{p,q}^{\left( r\right) }\left( r,z,\chi \right) =\underset%
{s\rightarrow r}{\lim }L_{p,q}^{\left( r\right) }\left( s,z,\chi \right) \\
&=&\underset{s\rightarrow r}{\lim }\frac{1}{\left[ F\right] _{q}^{r}}\sum_{%
\underset{\left( a_{1}+\cdots +a_{r},p\right) =1}{a_{1},\ldots ,a_{r}=1}%
}^{F}\chi \left( a_{1}+\cdots +a_{r}\right) \frac{\left\langle a_{1}+\cdots
+a_{r}+p^{\ast }z\right\rangle _{q}^{-s+r}}{\prod\limits_{j=1}^{r}\left(
s-j\right) } \\
&\times& \sum_{k=0}^{r-1}\binom{-s+r}{k}\beta _{k,q^{F}}^{\left(
r\right)
}q^{k\left( a_{1}+\cdots +a_{r}+p^{\ast }z\right) }\left( \frac{\left[ F%
\right] _{q}}{\left[ a_{1}+\cdots +a_{r}+p^{\ast }z\right]
_{q}}\right) ^{k}
\end{eqnarray*}
\begin{eqnarray*}
&&\hspace{-0.5in}+\underset{s\rightarrow r}{\lim }\frac{1}{\left[ F\right] _{q}^{r}}\sum_{%
\underset{\left( a_{1}+\cdots +a_{r},p\right) =1}{a_{1},\ldots ,a_{r}=1}%
}^{F}\chi \left( a_{1}+\cdots +a_{r}\right) \left\langle a_{1}+\cdots
+a_{r}+p^{\ast }z\right\rangle _{q}^{-s+r} \\
&&\times \sum_{k=r}^{\infty }\frac{\binom{-s+r}{k}}{\prod\limits_{j=1}^{r}%
\left( s-j\right) }\beta _{k,q^{F}}^{\left( r\right) }q^{k\left(
a_{1}+\cdots +a_{r}+p^{\ast }z\right) }\left( \frac{\left[ F\right] _{q}}{%
\left[ a_{1}+\cdots +a_{r}+p^{\ast }z\right] _{q}}\right) ^{k}. \\
&&\hspace{-0.5in}=\frac{1}{\left[ F\right] _{q}^{r}}\frac{1}{\left( r-1\right) !}\sum_{%
\underset{\left( a_{1}+\cdots +a_{r},p\right) =1}{a_{1},\ldots ,a_{r}=1}%
}^{F}\chi \left( a_{1}+\cdots +a_{r}\right) \left\{-\beta
_{0,q^{F}}^{\left( r\right) }\log _{p}\left\langle a_{1}+\cdots
+a_{r}+p^{\ast }z\right\rangle _{q}\right. \\
&&+\left. \sum_{k=r}^{\infty }\frac{\left( -1\right) ^{k}}{k}\beta
_{k,q^{F}}^{\left( r\right) }q^{k\left( a_{1}+\cdots
+a_{r}+p^{\ast }z\right) }\left( \frac{\left[ F\right]
_{q}}{\left[ a_{1}+\cdots +a_{r}+p^{\ast }z\right] _{q}}\right)
^{k}\right\} ,
\end{eqnarray*}

\noindent the desired result.
\end{proof}

\noindent Note that for $r=1$, Theorem \ref{thm4.1} reduces to Theorem 8 of
Kim \cite{Kim2005}.

\section{Congruences for Multiple Generalized $q$-Bernoulli Polynomials}

\hspace{0.15in}Congruences related to classical and generalized Bernoulli
numbers have found an amount of interest. One of the most celebrated
examples is the Kummer congruences for classical Bernoulli numbers (cf. \cite%
{Washington1997}):%
\begin{equation*}
p^{-1}\Delta _{c}\frac{B_{n}}{n}\in \mathbb{Z}_{p},
\end{equation*}

\noindent where $c\in \mathbb{Z}$, $c\geqslant 1$, $c\equiv 0\left(
mod\left( p-1\right) \right) $, and $n\in \mathbb{Z}$ is positive, even and $%
n\not\equiv 0\left( mod\left( p-1\right) \right) $. Here, $\Delta _{c}$ is
the forward difference operator which operates on a sequence $\left\{
x_{n}\right\} $ by%
\begin{equation*}
\Delta _{c}x_{n}=x_{n+c}-x_{n}.
\end{equation*}%
\noindent The powers $\Delta _{c}^{k}$ of $\Delta _{c}$ are defined by $%
\Delta _{c}^{0}=$identity and $\Delta _{c}^{k}=\Delta _{c}\circ \Delta
_{c}^{k-1}$ for positive integers $k$, so that%
\begin{equation*}
\Delta _{c}^{k}x_{n}=\sum_{m=0}^{k}\binom{k}{m}\left( -1\right)
^{k-m}x_{n+mc}.
\end{equation*}

\noindent More generally, it can be shown that%
\begin{equation*}
p^{-k}\Delta _{c}^{k}\frac{B_{n}}{n}\in \mathbb{Z}_{p},
\end{equation*}

\noindent where $k\in \mathbb{Z}$, $k\geqslant 1$ and $c$ and $n$ are as
above, but with $n>k$.

Kummer congruences for generalized Bernoulli numbers $B_{n,\chi }$ was first
regarded by Carlitz \cite{Carlitz1959}:

For positive $c\in \mathbb{Z}$, $c\equiv 0\left( mod\left( p-1\right)
\right) $, $n,k\in \mathbb{Z}$, $n>k\geqslant 1$, and $\chi $ such that $%
f=f_{\chi }\neq p^{m}$, where $m\in \mathbb{Z}$, $m\geqslant 0$,%
\begin{equation*}
p^{-k}\Delta _{c}^{k}\frac{B_{n,\chi }}{n}\in \mathbb{Z}_{p}\left[ \chi %
\right] .
\end{equation*}

\noindent Here, $\mathbb{Z}_{p}\left[ \chi \right] $ denotes the ring of
polynomials in $\chi $, whose coefficients are in $\mathbb{Z}_{p}$.

Shiratani \cite{Shiratani} applied the operator $\Delta _{c}^{k}$ to $%
-\left( 1-\chi _{n}\left( p\right) p^{n-1}\right) B_{n,\chi _{n}}/n$ for
similar $c$ and $\chi $, and showed that Carlitz's congruence is still true
without the restriction $n>k$, requiring only that $n\geqslant 1$. He also
established that the divisibility conditions on $c$ can be removed, and
proved%
\begin{equation*}
\left( p^{\ast }\right) ^{-k}\Delta _{c}^{k}\left( 1-\chi _{n}\left(
p\right) p^{n-1}\right) \frac{B_{n,\chi _{n}}}{n}\in \mathbb{Z}_{p}\left[
\chi \right] .
\end{equation*}

As an extension of the Kummer congruence, Gunaratne \cite%
{Gunaratne1995a,Gunaratne1995b} showed that the value%
\begin{equation*}
p^{-k}\Delta _{c}^{k}\left( 1-\chi _{n}\left( p\right) p^{n-1}\right) \frac{%
B_{n,\chi _{n}}}{n},
\end{equation*}

\noindent modulo $p\mathbb{Z}_{p}$, is independent of $n$ and%
\begin{equation*}
p^{-k}\Delta _{c}^{k}\left( 1-\chi _{n}\left( p\right) p^{n-1}\right) \frac{%
B_{n,\chi _{n}}}{n}\equiv p^{-k^{\prime }}\Delta _{c}^{k^{\prime }}\left(
1-\chi _{n^{\prime }}\left( p\right) p^{n^{\prime }-1}\right) \frac{%
B_{n^{\prime },\chi _{n^{\prime }}}}{n^{\prime }}\left( modp\mathbb{Z}%
_{p}\right) ,
\end{equation*}

\noindent if $p>3$, $c,n,k\in \mathbb{Z}$ are positive, $\chi =\omega ^{h}$,
where $h\in \mathbb{Z}$, $h\not\equiv 0\left( mod\left( p-1\right) \right) $%
, $n^{\prime },k^{\prime }\in \mathbb{Z}$, $k\equiv k^{\prime }\left(
mod\left( p-1\right) \right) $. Furthermore, by means of the binomial
coefficient operator%
\begin{equation*}
\binom{p^{-1}\Delta _{c}}{k}x_{n}=\frac{1}{k!}\left(
\prod\limits_{j=0}^{k-1}\left( p^{-1}\Delta _{c}-j\right) \right) x_{n},
\end{equation*}

\noindent it has been shown that for similar character $\chi $,%
\begin{equation*}
\binom{p^{-1}\Delta _{c}}{k}\left( 1-\chi _{n}\left( p\right) p^{n-1}\right)
\frac{B_{n,\chi _{n}}}{n}\in \mathbb{Z}_{p},
\end{equation*}

\noindent and this value, modulo $p\mathbb{Z}_{p}$, is independent of $n$.

Fox \cite{Fox2000} derived congruences similar to those above for the
generalized Bernoulli polynomials without restrictions on the character $%
\chi $. In \cite{Kim2002c}, Kim gave a proof of Kummer-type
congruence for the $q$-Bernoulli numbers of higher order. In
\cite{Cenkci-Kurt}, Cenkci and Kurt extended Fox's and Kim's
results to generalized $q$-Bernoulli polynomials. For other
versions of Kummer-type congruences related other numbers and
polynomials we refer \cite{Cenkci-Can-Kurt},
\cite{Simsek-Kim-Park-Ro-Jang-Rim}.

In order to derive a collection of congruences, similar to the results
above, relating to the multiple generalized $q$-Bernoulli polynomials, we
utilize the following theorem and its immediate consequence, found in \cite%
{Cenkci-Simsek-Kurt}:

\begin{theorem}
\label{thm5.1}For the character $\chi $ of conductor $f$, let $%
F_{0}=lcm\left( f,p^{\ast }\right) $, $F$ be a positive integer multiple of $%
p\left( p^{\ast }\right) ^{-1}rF_{0}$, $r\in \mathbb{Z}$, $r\geqslant 1$, $%
z\in \mathbb{C}_{p}$, $\left| z\right| _{p}\leqslant 1$ and $s\in D$, except
$s\neq 1,2,\ldots ,r$ if $\chi =1$. Then%
\begin{eqnarray*}
&&\hspace{-1in}L_{p,q}^{\left( r\right) }\left( s,z+rF,\chi \right)
-L_{p,q}^{\left( r\right) }\left( s,z,\chi \right) \\
&=&-\sum_{\underset{\left( a_{1}+\cdots +a_{r},p\right) =1}{a_{1},\ldots
,a_{r}=1}}^{p^{\ast }F}\chi _{r}\left( a_{1}+\cdots +a_{r}\right)
q^{a_{1}+\cdots +a_{r}+p^{\ast }z}\left\langle a_{1}+\cdots +a_{r}+p^{\ast
}z\right\rangle _{q}^{-s}.
\end{eqnarray*}
\end{theorem}

\begin{corollary}
\label{cor5.2}Let $F$, $r$ and $s$ be as in Theorem \ref{thm5.1}. Then%
\begin{eqnarray*}
&&\hspace{-1in}L_{p,q}^{\left( r\right) }\left( s,rF,\chi \right)
-L_{p,q}^{\left( r\right) }\left( s,\chi \right) \\
&=&-\sum_{\underset{\left( a_{1}+\cdots +a_{r},p\right) =1}{a_{1},\ldots
,a_{r}=1}}^{p^{\ast }F}\chi _{r}\left( a_{1}+\cdots +a_{r}\right)
q^{a_{1}+\cdots +a_{r}}\left\langle a_{1}+\cdots +a_{r}\right\rangle
_{q}^{-s},
\end{eqnarray*}

\noindent where $L_{p,q}^{\left( r\right) }\left( s,\chi \right)
=L_{p,q}^{\left( r\right) }\left( s,0,\chi \right) $.
\end{corollary}

\noindent We also incorporate the polynomial structure%
\begin{equation*}
B_{n}^{r}\left( z,q,\chi \right) =\frac{\left( -1\right) ^{r}}{\left(
n+1\right) _{r}}\left( \beta _{n+r,q,\chi _{n+r}}^{\left( r\right) }\left(
p^{\ast }z\right) -\chi _{n+r}\left( p\right) \left[ p\right] _{q}^{n}\beta
_{n+r,q^{p},\chi _{n+r}}^{\left( r\right) }\left( p^{-1}p^{\ast }z\right)
\right)
\end{equation*}%
\noindent and the set structure%
\begin{equation*}
R^{\ast }=\left\{ x\in \mathbb{Z}_{p}:\left| x\right| _{p}<p^{-1/\left(
p-1\right) }\right\} .
\end{equation*}%
\noindent Throughout, we assume that $q\in \mathbb{Z}_{p}$ with $\left|
1-q\right| _{p}<p^{-1/\left( p-1\right) }$, so that $q\equiv 1\left(
modR^{\ast }\right) $.

\begin{theorem}
\label{thm5.3}Let $n$, $c$, $k$, $r$ be positive integers and $z\in p\left(
p^{\ast }\right) ^{-1}rF_{0}R^{\ast }$. Then the quantity%
\begin{equation*}
\left( p^{\ast }\right) ^{-k}\Delta _{c}^{k}B_{n}^{r}\left( z,q,\chi \right)
-\left( p^{\ast }\right) ^{-k}\Delta _{c}^{k}B_{n}^{r}\left( 0,q,\chi
\right) \in R^{\ast }\left[ \chi \right] ,
\end{equation*}%
\noindent and, modulo $p^{\ast }R^{\ast }\left[ \chi \right] $, is
independent of $n$.
\end{theorem}

\begin{proof}
Since $\Delta _{c}$ is a linear operator, Corollary \ref{cor5.2} implies that%
\begin{eqnarray*}
&&\hspace{-0.5in}\Delta _{c}^{k}L_{p,q}^{\left( r\right) }\left( -n,rF,\chi
\right) -\Delta _{c}^{k}L_{p,q}^{\left( r\right) }\left( -n,\chi \right) \\
&=&-\sum_{\underset{\left( a_{1}+\cdots +a_{r},p\right) =1}{a_{1},\ldots
,a_{r}=1}}^{p^{\ast }F}\chi _{r}\left( a_{1}+\cdots +a_{r}\right)
q^{a_{1}+\cdots +a_{r}}\Delta _{c}^{k}\left\langle a_{1}+\cdots
+a_{r}\right\rangle _{q}^{n}.
\end{eqnarray*}%
\noindent Thus%
\begin{eqnarray*}
&&\hspace{-0.5in}\Delta _{c}^{k}B_{n}^{r}\left( rF,q,\chi \right) -\Delta
_{c}^{k}B_{n}^{r}\left( 0,q,\chi \right) \\
&=&-\sum_{\underset{\left( a_{1}+\cdots +a_{r},p\right) =1}{a_{1},\ldots
,a_{r}=1}}^{p^{\ast }F}\chi _{r}\left( a_{1}+\cdots +a_{r}\right)
q^{a_{1}+\cdots +a_{r}}\Delta _{c}^{k}\left\langle a_{1}+\cdots
+a_{r}\right\rangle _{q}^{n}.
\end{eqnarray*}%
\noindent Note that%
\begin{eqnarray*}
\Delta _{c}^{k}\left\langle a_{1}+\cdots +a_{r}\right\rangle _{q}^{n}
&=&\sum_{m=0}^{k}\binom{k}{m}\left( -1\right) ^{k-m}\left\langle
a_{1}+\cdots +a_{r}\right\rangle _{q}^{n+mc} \\
&=&\left\langle a_{1}+\cdots +a_{r}\right\rangle _{q}^{n}\left( \left\langle
a_{1}+\cdots +a_{r}\right\rangle _{q}^{c}-1\right) ^{k}.
\end{eqnarray*}%
\noindent Now, $\left\langle a_{1}+\cdots +a_{r}\right\rangle _{q}\equiv
1\left( modp^{\ast }R^{\ast }\right) $ for $\left( a_{1}+\cdots
+a_{r},p\right) =1$, which implies that $\left\langle a_{1}+\cdots
+a_{r}\right\rangle _{q}^{c}\equiv 1\left( modp^{\ast }R^{\ast }\right) $,
and thus $\Delta _{c}^{k}\left\langle a_{1}+\cdots +a_{r}\right\rangle
_{q}^{n}\equiv 0\left( mod\left( p^{\ast }\right) ^{k}R^{\ast }\right) $.
Therefore%
\begin{equation*}
\Delta _{c}^{k}B_{n}^{r}\left( rF,q,\chi \right) -\Delta
_{c}^{k}B_{n}^{r}\left( 0,q,\chi \right) \equiv 0\left( mod\left( p^{\ast
}\right) ^{k}R^{\ast }\left[ \chi \right] \right) ,
\end{equation*}%
\noindent and so%
\begin{equation*}
\left( p^{\ast }\right) ^{-k}\Delta _{c}^{k}B_{n}^{r}\left( rF,q,\chi
\right) -\left( p^{\ast }\right) ^{-k}\Delta _{c}^{k}B_{n}^{r}\left(
0,q,\chi \right) \in R^{\ast }\left[ \chi \right] .
\end{equation*}%
\noindent Also, since $\left\langle a_{1}+\cdots +a_{r}\right\rangle
_{q}^{c}\equiv 1\left( modp^{\ast }R^{\ast }\right) $, the equation%
\begin{eqnarray*}
&&\hspace{-0.5in}\Delta _{c}^{k}B_{n}^{r}\left( rF,q,\chi \right) -\Delta
_{c}^{k}B_{n}^{r}\left( 0,q,\chi \right) \\
&=&-\sum_{\underset{\left( a_{1}+\cdots +a_{r},p\right) =1}{a_{1},\ldots
,a_{r}=1}}^{p^{\ast }F}\chi _{r}\left( a_{1}+\cdots +a_{r}\right)
q^{a_{1}+\cdots +a_{r}}\left\langle a_{1}+\cdots +a_{r}\right\rangle
_{q}^{n}\left( \frac{\left\langle a_{1}+\cdots +a_{r}\right\rangle _{q}^{c}-1%
}{p^{\ast }}\right) ^{k}
\end{eqnarray*}%
\noindent implies that the value of $\left( p^{\ast }\right) ^{-k}\Delta
_{c}^{k}B_{n}^{r}\left( rF,q,\chi \right) -\left( p^{\ast }\right)
^{-k}\Delta _{c}^{k}B_{n}^{r}\left( 0,q,\chi \right) $, modulo $p^{\ast
}R^{\ast }\left[ \chi \right] $, is independent of $n$.

Let $z\in p\left( p^{\ast }\right) ^{-1}rF_{0}R^{\ast }$. Since the set of
positive integers in $p\left( p^{\ast }\right) ^{-1}rF_{0}\mathbb{Z}$ is
dense in \linebreak$p\left( p^{\ast }\right) ^{-1}rF_{0}R^{\ast }$, there
exists a sequence $\left\{ z_{j}\right\} $ in $p\left( p^{\ast }\right)
^{-1}rF_{0}\mathbb{Z}$ with $z_{j}>0$ for each $j$, such that $%
z_{j}\rightarrow z$. Now, $B_{n}^{r}\left( z,q,\chi \right) $ is a
polynomial, which implies that $B_{n}^{r}\left( z_{j},q,\chi \right)
\rightarrow B_{n}^{r}\left( z,q,\chi \right) $. Therefore,%
\begin{equation*}
\underset{j\rightarrow \infty }{\lim }\left( \Delta _{c}^{k}B_{n}^{r}\left(
z_{j},q,\chi \right) -\Delta _{c}^{k}B_{n}^{r}\left( 0,q,\chi \right)
\right) =\Delta _{c}^{k}B_{n}^{r}\left( z,q,\chi \right) -\Delta
_{c}^{k}B_{n}^{r}\left( 0,q,\chi \right) .
\end{equation*}%
\noindent The left side of this equality is $0$ modulo $\left( p^{\ast
}\right) ^{k}R^{\ast }\left[ \chi \right] $, which implies that%
\begin{equation*}
\Delta _{c}^{k}B_{n}^{r}\left( z,q,\chi \right) -\Delta
_{c}^{k}B_{n}^{r}\left( 0,q,\chi \right) \equiv 0\left( mod\left( p^{\ast
}\right) ^{k}R^{\ast }\left[ \chi \right] \right) ,
\end{equation*}%
\noindent and so%
\begin{equation*}
\left( p^{\ast }\right) ^{-k}\Delta _{c}^{k}B_{n}^{r}\left( z,q,\chi \right)
-\left( p^{\ast }\right) ^{-k}\Delta _{c}^{k}B_{n}^{r}\left( 0,q,\chi
\right) \in R^{\ast }\left[ \chi \right] .
\end{equation*}%
\noindent Furthermore, for a positive integer $n^{\prime }$%
\begin{eqnarray*}
&&\hspace{-0.5in}\underset{j\rightarrow \infty }{\lim }\left\{ \left( \left(
p^{\ast }\right) ^{-k}\Delta _{c}^{k}B_{n}^{r}\left( z_{j},q,\chi \right)
-\left( p^{\ast }\right) ^{-k}\Delta _{c}^{k}B_{n}^{r}\left( 0,q,\chi
\right) \right) \right. \\
&&\left. -\left( \left( p^{\ast }\right) ^{-k}\Delta _{c}^{k}B_{n^{\prime
}}^{r}\left( z_{j},q,\chi \right) -\left( p^{\ast }\right) ^{-k}\Delta
_{c}^{k}B_{n^{\prime }}^{r}\left( 0,q,\chi \right) \right) \right\} \\
&&\hspace{-0.5in}=\left\{ \left( \left( p^{\ast }\right) ^{-k}\Delta
_{c}^{k}B_{n}^{r}\left( z,q,\chi \right) -\left( p^{\ast }\right)
^{-k}\Delta _{c}^{k}B_{n}^{r}\left( 0,q,\chi \right) \right) \right. \\
&&\left. -\left( \left( p^{\ast }\right) ^{-k}\Delta _{c}^{k}B_{n^{\prime
}}^{r}\left( z,q,\chi \right) -\left( p^{\ast }\right) ^{-k}\Delta
_{c}^{k}B_{n^{\prime }}^{r}\left( 0,q,\chi \right) \right) \right\} .
\end{eqnarray*}%
\noindent Since $z_{j}\in p\left( p^{\ast }\right) ^{-1}rF_{0}\mathbb{Z}$
for all $j$, the quantity on the left must be $0$ modulo $p^{\ast }R^{\ast }%
\left[ \chi \right] $. Therefore, the value $\left( p^{\ast }\right)
^{-k}\Delta _{c}^{k}B_{n}^{r}\left( z,q,\chi \right) -\left( p^{\ast
}\right) ^{-k}\Delta _{c}^{k}B_{n}^{r}\left( 0,q,\chi \right) $, modulo $%
p^{\ast }R^{\ast }\left[ \chi \right] $, is independent of $n$.
\end{proof}

\begin{theorem}
\label{thm5.4}Let $n$, $c$, $r$, $k$, $k^{\prime }$ be positive integers, $%
k\equiv k^{\prime }\left( mod\left( p-1\right) \right) $ and let $z\in
p\left( p^{\ast }\right) ^{-1}rF_{0}R^{\ast }$. Then%
\begin{eqnarray*}
&&\hspace{-0.5in}\left( p^{\ast }\right) ^{-k}\Delta _{c}^{k}B_{n}^{r}\left(
z,q,\chi \right) -\left( p^{\ast }\right) ^{-k}\Delta
_{c}^{k}B_{n}^{r}\left( 0,q,\chi \right) \\
&\equiv &\left( p^{\ast }\right) ^{-k^{\prime }}\Delta _{c}^{k^{\prime
}}B_{n}^{r}\left( z,q,\chi \right) -\left( p^{\ast }\right) ^{-k^{\prime
}}\Delta _{c}^{k^{\prime }}B_{n}^{r}\left( 0,q,\chi \right) \left(
modpR^{\ast }\left[ \chi \right] \right) .
\end{eqnarray*}
\end{theorem}

\begin{proof}
Let $k$ and $k^{\prime }$ be positive integers such that $k\equiv k^{\prime
}\left( mod\left( p-1\right) \right) $. Without loss of generality, assume
that $k\geqslant k^{\prime }$. Then%
\begin{eqnarray*}
&&\hspace{-1in}\left( \left( p^{\ast }\right) ^{-k}\Delta
_{c}^{k}B_{n}^{r}\left( rF,q,\chi \right) -\left( p^{\ast }\right)
^{-k}\Delta _{c}^{k}B_{n}^{r}\left( 0,q,\chi \right) \right)  \\
&-&\left( \left( p^{\ast }\right) ^{-k^{\prime }}\Delta _{c}^{k^{\prime
}}B_{n}^{r}\left( rF,q,\chi \right) -\left( p^{\ast }\right) ^{-k^{\prime
}}\Delta _{c}^{k^{\prime }}B_{n}^{r}\left( 0,q,\chi \right) \right)
\end{eqnarray*}%
\begin{eqnarray*}
&&\hspace{-0.3in}=-\sum_{\underset{\left( a_{1}+\cdots +a_{r},p\right) =1}{%
a_{1},\ldots ,a_{r}=1}}^{p^{\ast }F}\chi _{r}\left( a_{1}+\cdots
+a_{r}\right) q^{a_{1}+\cdots +a_{r}}\left\langle a_{1}+\cdots
+a_{r}\right\rangle _{q}^{n} \\
&&\hspace{-0.2in}\times \left\{ \left( \frac{\left\langle a_{1}+\cdots
+a_{r}\right\rangle _{q}^{c}-1}{p^{\ast }}\right) ^{k}-\left( \frac{%
\left\langle a_{1}+\cdots +a_{r}\right\rangle _{q}^{c}-1}{p^{\ast }}\right)
^{k^{\prime }}\right\}  \\
&&\hspace{-0.3in}=-\sum_{\underset{\left( a_{1}+\cdots +a_{r},p\right) =1}{%
a_{1},\ldots ,a_{r}=1}}^{p^{\ast }F}\chi _{r}\left( a_{1}+\cdots
+a_{r}\right) q^{a_{1}+\cdots +a_{r}}\left\langle a_{1}+\cdots
+a_{r}\right\rangle _{q}^{n} \\
&&\hspace{-0.2in}\times \left( \frac{\left\langle a_{1}+\cdots
+a_{r}\right\rangle _{q}^{c}-1}{p^{\ast }}\right) ^{k^{\prime }}\left\{
\left( \frac{\left\langle a_{1}+\cdots +a_{r}\right\rangle _{q}^{c}-1}{%
p^{\ast }}\right) ^{k-k^{\prime }}-1\right\} .
\end{eqnarray*}%
\noindent If $a_{1}+\cdots +a_{r}$ such that%
\begin{equation*}
\left\langle a_{1}+\cdots +a_{r}\right\rangle _{q}^{c}-1\not\equiv 0\left(
modpp^{\ast }R^{\ast }\right) ,
\end{equation*}%
\noindent then, since $k-k^{\prime }\equiv 0\left( \text{mod}\left(
p-1\right) \right) $, we have%
\begin{equation*}
\left( \frac{\left\langle a_{1}+\cdots +a_{r}\right\rangle _{q}^{c}-1}{%
p^{\ast }}\right) ^{k-k^{\prime }}-1\equiv 0\left( modpR^{\ast }\right) .
\end{equation*}%
\noindent Thus%
\begin{eqnarray*}
&&\hspace{-0.5in}\left( p^{\ast }\right) ^{-k}\Delta _{c}^{k}B_{n}^{r}\left(
rF,q,\chi \right) -\left( p^{\ast }\right) ^{-k}\Delta
_{c}^{k}B_{n}^{r}\left( 0,q,\chi \right)  \\
&\equiv &\left( p^{\ast }\right) ^{-k^{\prime }}\Delta _{c}^{k^{\prime
}}B_{n}^{r}\left( rF,q,\chi \right) -\left( p^{\ast }\right) ^{-k^{\prime
}}\Delta _{c}^{k^{\prime }}B_{n}^{r}\left( 0,q,\chi \right) \left(
modpR^{\ast }\left[ \chi \right] \right) .
\end{eqnarray*}

Now let $z\in p\left( p^{\ast }\right) ^{-1}rF_{0}R^{\ast }$. Then there
exists a sequence $\left\{ z_{j}\right\} $ in $p\left( p^{\ast }\right)
^{-1}rF_{0}\mathbb{Z}$ with $z_{j}>0$ for each $j$, such that $%
z_{j}\rightarrow z$. Consider%
\begin{eqnarray*}
&&\hspace{-0.5in}\underset{j\rightarrow \infty }{\lim }\left\{ \left( \left(
p^{\ast }\right) ^{-k}\Delta _{c}^{k}B_{n}^{r}\left( z_{j},q,\chi \right)
-\left( p^{\ast }\right) ^{-k}\Delta _{c}^{k}B_{n}^{r}\left( 0,q,\chi
\right) \right) \right. \\
&&\left. -\left( \left( p^{\ast }\right) ^{-k^{\prime }}\Delta
_{c}^{k^{\prime }}B_{n}^{r}\left( z_{j},q,\chi \right) -\left( p^{\ast
}\right) ^{-k^{\prime }}\Delta _{c}^{k^{\prime }}B_{n}^{r}\left( 0,q,\chi
\right) \right) \right\} \\
&&\hspace{-0.5in}=\left\{ \left( \left( p^{\ast }\right) ^{-k}\Delta
_{c}^{k}B_{n}^{r}\left( z,q,\chi \right) -\left( p^{\ast }\right)
^{-k}\Delta _{c}^{k}B_{n}^{r}\left( 0,q,\chi \right) \right) \right. \\
&&\left. -\left( \left( p^{\ast }\right) ^{-k^{\prime }}\Delta
_{c}^{k^{\prime }}B_{n}^{r}\left( z,q,\chi \right) -\left( p^{\ast }\right)
^{-k^{\prime }}\Delta _{c}^{k^{\prime }}B_{n}^{r}\left( 0,q,\chi \right)
\right) \right\} .
\end{eqnarray*}%
\noindent Since the left side of this equality must be $0$ modulo $pR^{\ast }%
\left[ \chi \right] $, the proof follows.
\end{proof}

The binomial coefficients operator $\binom{T}{k}$ associated to an operator $%
T$ is defined by writing the binomial coefficients%
\begin{equation*}
\binom{X}{k}=\frac{X\left( X-1\right) \cdots \left( X-k+1\right) }{k!}
\end{equation*}
\noindent for $k\geqslant 0$ as a polynomial in $X$, and replacing $X$ by $T$%
.

In the proof of next theorem, we need special numbers, namely the Stirling
numbers of the first kind $s\left( n,k\right) $, which are defined by means
of the generating function%
\begin{equation*}
\frac{\left( \text{log}\left( 1+t\right) \right) ^{k}}{k!}%
=\sum_{n=0}^{\infty }s\left( n,k\right) \frac{t^{n}}{n!},
\end{equation*}

\noindent for $k\in \mathbb{Z}$, $k\geqslant 0$. Since there is no constant
term in the expansion of log$\left( 1+t\right) $, $s\left( n,k\right) =0$
for $0\leqslant n<k$. Also, $s\left( n,n\right) =1$ for all $n\geqslant 0$.
The numbers $s\left( n,k\right) $ are integers and satisfy the following
relation related to binomial coefficients:%
\begin{equation}
\binom{x}{k}=\frac{1}{n!}\sum_{k=0}^{n}s\left( n,k\right) x^{k}.  \label{5.1}
\end{equation}

\noindent For further information for Stirling numbers, we refer \cite%
{Comtet}.

\begin{theorem}
\label{thmEk5}Let $n$, $c$, $k$ be positive integers and $z\in p\left(
p^{\ast }\right) ^{-1}F_{0}R^{\ast }$. Then the quantity%
\begin{equation*}
\binom{\left( p^{\ast }\right) ^{-1}\Delta _{c}}{k}B_{n}^{r}\left( z,q,\chi
\right) -\binom{\left( p^{\ast }\right) ^{-1}\Delta _{c}}{k}B_{n}^{r}\left(
0,q,\chi \right) \in R^{\ast }\left[ \chi \right] ,
\end{equation*}%
\noindent and, modulo $p^{\ast }R^{\ast }\left[ \chi \right] $, is
independent of $n$.
\end{theorem}

\begin{proof}
Since the binomial coefficients operator is a linear operator, Corollary \ref%
{cor5.2} implies that%
\begin{eqnarray*}
&&\hspace{-0.5in}\binom{\left( p^{\ast }\right) ^{-1}\Delta _{c}}{k}%
L_{p,q}^{\left( r\right) }\left( -n,rF,\chi \right) -\binom{\left( p^{\ast
}\right) ^{-1}\Delta _{c}}{k}L_{p,q}^{\left( r\right) }\left( -n,\chi \right)
\\
&=&-\sum_{\underset{\left( a_{1}+\cdots +a_{r},p\right) =1}{a_{1},\ldots
,a_{r}=1}}^{p^{\ast }F}\chi _{r}\left( a_{1}+\cdots +a_{r}\right)
q^{a_{1}+\cdots +a_{r}}\binom{\left( p^{\ast }\right) ^{-1}\Delta _{c}}{k}%
\left\langle a_{1}+\cdots +a_{r}\right\rangle _{q}^{n}.
\end{eqnarray*}%
\noindent Then%
\begin{eqnarray*}
&&\hspace{-0.5in}\binom{\left( p^{\ast }\right) ^{-1}\Delta _{c}}{k}%
B_{n}^{r}\left( rF,q,\chi \right) -\binom{\left( p^{\ast }\right)
^{-1}\Delta _{c}}{k}B_{n}^{r}\left( 0,q,\chi \right) \\
&=&-\sum_{\underset{\left( a_{1}+\cdots +a_{r},p\right) =1}{a_{1},\ldots
,a_{r}=1}}^{p^{\ast }F}\chi _{r}\left( a_{1}+\cdots +a_{r}\right)
q^{a_{1}+\cdots +a_{r}}\binom{\left( p^{\ast }\right) ^{-1}\Delta _{c}}{k}%
\left\langle a_{1}+\cdots +a_{r}\right\rangle _{q}^{n}.
\end{eqnarray*}%
\noindent By using (\ref{5.1}), we can write%
\begin{eqnarray*}
&&\hspace{-0.5in}\binom{\left( p^{\ast }\right) ^{-1}\Delta _{c}}{k}%
\left\langle a_{1}+\cdots +a_{r}\right\rangle _{q}^{n} \\
&=&\frac{1}{k!}\sum_{m=0}^{k}s\left( k,m\right) \left( p^{\ast }\right)
^{-m}\Delta _{c}^{m}\left\langle a_{1}+\cdots +a_{r}\right\rangle _{q}^{n} \\
&=&\frac{1}{k!}\sum_{m=0}^{k}s\left( k,m\right) \left( p^{\ast }\right)
^{-m}\left\langle a_{1}+\cdots +a_{r}\right\rangle _{q}^{n}\left(
\left\langle a_{1}+\cdots +a_{r}\right\rangle _{q}^{c}-1\right) ^{m}.
\end{eqnarray*}%
\noindent Thus%
\begin{eqnarray*}
&&\hspace{-1in}\binom{\left( p^{\ast }\right) ^{-1}\Delta _{c}}{k}%
B_{n}^{r}\left( rF,q,\chi \right) -\binom{\left( p^{\ast }\right)
^{-1}\Delta _{c}}{k}B_{n}^{r}\left( 0,q,\chi \right) \\
&=&-\sum_{\underset{\left( a_{1}+\cdots +a_{r},p\right) =1}{a_{1},\ldots
,a_{r}=1}}^{p^{\ast }F}\chi _{r}\left( a_{1}+\cdots +a_{r}\right)
q^{a_{1}+\cdots +a_{r}} \\
&&\hspace{-0.5in}\times \binom{\left( p^{\ast }\right) ^{-1}\Delta _{c}}{k}%
\left\langle a_{1}+\cdots +a_{r}\right\rangle _{q}^{n}\binom{\left( p^{\ast
}\right) ^{-1}\left( \left\langle a_{1}+\cdots +a_{r}\right\rangle
_{q}^{c}-1\right) }{k}.
\end{eqnarray*}%
\noindent Since $\left( p^{\ast }\right) ^{-1}\left( \left\langle
a_{1}+\cdots +a_{r}\right\rangle _{q}^{c}-1\right) \in R^{\ast }$ for $%
\left( a_{1}+\cdots +a_{r},p\right) =1$, we see that%
\begin{equation*}
\left\langle a_{1}+\cdots +a_{r}\right\rangle _{q}^{n}\binom{\left( p^{\ast
}\right) ^{-1}\left( \left\langle a_{1}+\cdots +a_{r}\right\rangle
_{q}^{c}-1\right) }{k}\in R^{\ast }.
\end{equation*}%
\noindent This then implies%
\begin{equation*}
\binom{\left( p^{\ast }\right) ^{-1}\Delta _{c}}{k}B_{n}^{r}\left( rF,q,\chi
\right) -\binom{\left( p^{\ast }\right) ^{-1}\Delta _{c}}{k}B_{n}^{r}\left(
0,q,\chi \right) \in R^{\ast }\left[ \chi \right] .
\end{equation*}%
\noindent Furthermore, since $\left\langle a_{1}+\cdots +a_{r}\right\rangle
_{q}^{c}\equiv 1\left( modp^{\ast }R^{\ast }\right) $, the value of this
quantity, modulo $p^{\ast }R^{\ast }\left[ \chi \right] $, is independent of
$n$.

Now let $z\in p\left( p^{\ast }\right) ^{-1}rF_{0}R^{\ast }$, and let $%
\left\{ z_{j}\right\} $ be a sequence in $p\left( p^{\ast }\right)
^{-1}rF_{0}\mathbb{Z}$, with $z_{j}>0$ for each $j$, such that $%
z_{j}\rightarrow z$. Then%
\begin{eqnarray*}
&&\hspace{-1.5in}\underset{j\rightarrow \infty }{\lim }\binom{\left( p^{\ast
}\right) ^{-1}\Delta _{c}}{k}B_{n}^{r}\left( z_{j},q,\chi \right) -\binom{%
\left( p^{\ast }\right) ^{-1}\Delta _{c}}{k}B_{n}^{r}\left( 0,q,\chi \right)
\\
&=&\binom{\left( p^{\ast }\right) ^{-1}\Delta _{c}}{k}B_{n}^{r}\left(
z,q,\chi \right) -\binom{\left( p^{\ast }\right) ^{-1}\Delta _{c}}{k}%
B_{n}^{r}\left( 0,q,\chi \right)
\end{eqnarray*}%
\noindent must be in $R^{\ast }\left[ \chi \right] $. Now let $n^{\prime
}\in \mathbb{Z}$, $n^{\prime }>0$, and consider

\begin{eqnarray*}
&&\hspace{-0.5in}\underset{j\rightarrow \infty }{\text{lim}}\left\{ \binom{%
\left( p^{\ast }\right) ^{-1}\Delta _{c}}{k}B_{n}^{r}\left( z_{j},q,\chi
\right) -\binom{\left( p^{\ast }\right) ^{-1}\Delta _{c}}{k}B_{n}^{r}\left(
0,q,\chi \right) \right. \\
&&\left. -\binom{\left( p^{\ast }\right) ^{-1}\Delta _{c}}{k}B_{n^{\prime
}}^{r}\left( z_{j},q,\chi \right) -\binom{\left( p^{\ast }\right)
^{-1}\Delta _{c}}{k}B_{n^{\prime }}^{r}\left( 0,q,\chi \right) \right\} \\
&&\hspace{-0.5in}=\left\{ \binom{\left( p^{\ast }\right) ^{-1}\Delta _{c}}{k}%
B_{n}^{r}\left( z,q,\chi \right) -\binom{\left( p^{\ast }\right) ^{-1}\Delta
_{c}}{k}B_{n}^{r}\left( 0,q,\chi \right) \right. \\
&&\left. -\binom{\left( p^{\ast }\right) ^{-1}\Delta _{c}}{k}B_{n^{\prime
}}^{r}\left( z,q,\chi \right) -\binom{\left( p^{\ast }\right) ^{-1}\Delta
_{c}}{k}B_{n^{\prime }}^{r}\left( 0,q,\chi \right) \right\} .
\end{eqnarray*}%
\noindent The quantity on the left must be $0$ modulo $p^{\ast }R^{\ast }%
\left[ \chi \right] $, which implies that the value of%
\begin{equation*}
\binom{\left( p^{\ast }\right) ^{-1}\Delta _{c}}{k}B_{n}^{r}\left( z,q,\chi
\right) -\binom{\left( p^{\ast }\right) ^{-1}\Delta _{c}}{k}B_{n}^{r}\left(
0,q,\chi \right) ,
\end{equation*}%
\noindent modulo $p^{\ast }R^{\ast }\left[ \chi \right] $, is independent of
$n$.
\end{proof}

\end{document}